\documentclass[a4paper,11pt]{article}
\usepackage{amsfonts}
\usepackage{mathrsfs}
\usepackage{amsmath,amssymb,amsthm}
\usepackage{latexsym}
\usepackage{indentfirst}
\usepackage[top=1in, bottom=1in, left=1in, right=1in]{geometry}
\usepackage{color}
\usepackage{footmisc}
\usepackage{endnotes}
\usepackage{algorithm}
\usepackage{float}
\usepackage{graphicx}
\usepackage{epstopdf}
\usepackage{longtable}
\usepackage{authblk}
 \usepackage[authoryear]{natbib}

\def\bxi{\boldsymbol{\Xi}}
\newcommand{\ud}{\mathrm{d}}

\newcommand{\DP}{\mathrm{DP}}
\newcommand{\NGGP}{\mathrm{NGGP}}

\def\bt{\boldsymbol{\theta}}
\def\bxi{\boldsymbol{\xi}}

\newtheorem{theorem}{Theorem}
\newtheorem{lemma}{Lemma}

\newtheorem{prop}{Proposition}
\newtheorem{definition}{Definition}
\theoremstyle{definition}
\newtheorem{remark}{Remark}

\allowdisplaybreaks

\begin{document}

\title{\textbf{On Posterior Consistency of Tail Index for Bayesian Kernel Mixture Models}}

\author[1]{Cheng Li \thanks{stalic@nus.edu.sg}}
\author[2]{Lizhen Lin \thanks{lizhen.lin@nd.edu}}
\author[3]{David B. Dunson \thanks{dunson@duke.edu}}

\affil[1]{Department of Statistics and Applied Probability, National University of Singapore}
\affil[2]{Department of Applied and Computational Mathematics and Statistics, The University of Notre Dame}
\affil[3]{Department of Statistical Science, Duke University}

\date{}
\maketitle

\begin{abstract}
Asymptotic theory of tail index estimation has been studied extensively in the frequentist literature on extreme values, but rarely in the Bayesian context. We investigate whether popular Bayesian kernel mixture models  are able to support heavy tailed distributions and consistently estimate the tail index. We show that posterior inconsistency in tail index is surprisingly common for both parametric and nonparametric mixture models. We then present a set of sufficient conditions under which posterior consistency in tail index can be achieved, and verify these conditions for Pareto mixture models under general mixing priors.
\end{abstract}


\section{Introduction}
Datasets from a variety of fields, such as environmental science, finance, industrial engineering, and telecommunications, demonstrate heavy tailed behavior that can substantially influence statistical inference and decision making.  It is of interest to develop estimation methods that can capture both the bulk of the data and the tails accurately.
Bayesian kernel mixture models  provide  a flexible framework for  density estimation with strong large sample guarantees. Some of the most popular models include finite mixtures (MFM, \citealt{RicGre97, GreRic01}), Dirichlet process mixtures (DPM, \citealt{Fer73, Lo84, Mac94, EscWes95, Nea00}), and mixtures with mixing measures given by normalized random measures with independent increments (NRMI, \citealt{Regetal03, Lijetal07, Jametal09, LijPru10, Baretal13, FavTeh13}). However, most of the existing literature on Bayesian asymptotics for density estimation, including results on posterior consistency and convergence rates, assumes that the true density has either a compact support or exponentially decaying tails (\citealt{Ghoetal99, Ghoetal00, GhoVaa07, KruRou10, Sheetal13}). A few exceptions such as \citet{Tok06} and \citet{WuGho08} have shown posterior consistency for some heavy tailed densities for specific kernel mixture models. There exist fundamental limitations and barriers in understanding  the tail behavior of kernel mixture models and their large sample properties, especially for models with nonparametric mixing priors.


The current paper investigates theory on the tail behavior of popular Bayesian kernel mixture models, assessing whether they are suitable for modeling heavy tailed distributions. We focus on studying the tails of univariate continuous densities and assume that the true density has polynomially decaying tails. Such power law behavior has been observed in many real data applications (see \citealt{Claetal09} for a review).  Denote $f_0$ and $F_0$ as the true density function and the true cumulative distribution function (cdf) on $\mathbb{R}$. Let $\overline F(x) = 1-F(x)$ for $x\in \mathbb{R}$ be the survival function of $F$. For sufficiently large $x$, a distribution $F$ with a polynomially decaying right tail can be described by the relation
\begin{align}\label{righttail}
& \overline F(x) = x^{-\alpha_+(F)}L_F(x),
\end{align}
where $\alpha_+(F)>0$ is the \emph{right tail index}, and $L_F$ is a {\it slowly varying function} that satisfies $\lim_{y\to +\infty} L_F(xy)/L_F(y)=1$ for any $x>0$.  In this paper, we will only consider a true distribution $F_0$ that satisfies the relation \eqref{righttail}. The decay rate in the right tail of $F_0$ can be characterized by the right tail index $\alpha_+(F_0)$, up to some slowly varying function $L_{F_0}$.  The left tail index can be defined similarly.  If $\alpha_+(F_0)\in (0,+\infty)$, then from extreme value theory, the distribution $F_0$ falls within the class of \emph{Fr\'echet maximum domain of attraction} (FMDA, \citealt{Beietal04}). Examples of distributions satisfying \eqref{righttail} with $\alpha_+(F_0)\in (0,+\infty)$ include the Pareto distribution, the Student's $t$ distribution, the $F$ distribution, the inverse gamma distribution, the log-gamma distribution, the Burr distribution, etc.

We study theoretical properties of the posterior distribution of the right tail index $\alpha_+(F)$ in a Bayesian framework. In particular, we consider the kernel mixture model:
\begin{align}\label{model1}
f(x) = \int k\left(x;\bt \right) \ud G(\bt),\qquad G \sim \pi(G;\bxi),
\end{align}
where $k(\cdot;\bt)$ is a univariate kernel function with parameter $\bt$ such that $\int k(x;\bt)\ud x=1$ for all $\bt$, $G$ is a mixing measure of $\bt$, and $\pi$ is the prior on $G$ with hyperparameters $\bxi$.  This model is quite general, covering the aforementioned MFM, DPM and NRMI mixture models as special cases.

We will answer two critical questions for understanding how model \eqref{model1} can handle heavy tailed densities: (i) what choices of kernels and priors for the mixing measure can generate density functions with \emph{tail indices varying} in a reasonable range, and (ii) under  what types of conditions can one guarantee that the tail indices from the  posterior are close to the tail index of the true distribution. The first question is related to whether the Bayesian kernel mixture model is capable of flexibly fitting heavy tailed distributions with different decay rates. The second question is on the frequentist asymptotic properties of Bayesian models estimating the tail index, requiring substantial extension of the scope of existing theory for Bayesian density estimation.

There is a rich literature on frequentist estimation of the tail index.  Most of the estimators are constructed from tail order statistics, such as the Hill's estimator (\citealt{Hil75, DeHRes80}), the Pickands' estimator (\citealt{Pic75}) and their variations.  The Hill's estimator is consistent (\citealt{Mas82}) and asymptotically normal with appropriate choices of the tail order statistics for certain nonparametric classes of distributions (\citealt{Hal82, HaeTeu85}). Minimax rates for the tail index have been obtained under different classes of distributions (\citealt{HalWel84, Dre98, Dre01, Nov14, CarKim15}), and they are attainable through adaptive estimators (\citealt{HalWel85, CarKim15, BouTho15}).

However, there is a lack of understanding of the properties of likelihood-based approaches.  The limited Bayesian literature has focused on a peak-over-threshold (POT) strategy, with the tail of the density over a high threshold $t$ assumed to follow a generalized Pareto distribution.   If $F$ belongs to FDMA with right tail index $\alpha_+(F)$, then as the threshold $t$ becomes large, the right excess distribution $\widetilde F(x) = F(t+x)/F(t)$ for $x>0$ converges in law to a generalized Pareto distribution with tail index $\alpha_+(F)$. Posterior sampling schemes have been discussed in, for example, \cite{Frietal02, Botetal03, SteTaw04, Dieetal05, Nasetal12, Wanetal12, Fuq15}.  The POT strategy can be viewed as artificial in choosing different models below and above the threshold, with the restriction of a parametric tail.  The kernel mixture model allows one to choose a single flexible model for all of the data including the tails. \cite{Tre08} argues in favor of such an approach in using DPMs of Pareto kernels.

The rest of the paper is organized as follows. In Section \ref{prelim}, we formally introduce the definition of tail index and the posterior consistency of tail index. In Section \ref{lsm} we show that in general, location-scale kernel mixture models cannot generate densities with varying tail indices, even if the kernel is heavy tailed. In particular, our results reveal that in many cases, the posterior distribution under the mixture model can only generate distributions with a \emph{singleton index}. In Section \ref{suff}, we provide general sufficient conditions for Bayesian posterior consistency of tail index. These conditions are then verified for the example of Pareto kernel mixtures in Section \ref{example}. Section \ref{disc} concludes with discussions. Technical proofs are included in the appendix and the supplementary material.


\section{Preliminaries for Tail Index in a Bayesian Framework}\label{prelim}
\subsection{Definition of Tail Index}
We first describe a notion of tail index for any distribution $F$ with density $f$ defined on $\mathbb{R}$. For $x\in \mathbb{R}$, we define its right and left tail indices as
\begin{align}\label{tidef}
& \alpha_+(F) = \liminf_{x\to +\infty} \frac{-\log \overline F(x) }{\log x} = \liminf_{x\to +\infty} \frac{-\log P_F(X>x) }{\log x}, \\
& \alpha_-(F) = \liminf_{x\to -\infty} \frac{-\log F(x) }{\log (-x)}=\liminf_{x\to -\infty} \frac{-\log P_F(X\leq x) }{\log (-x)},\nonumber
\end{align}
where $P_F(\cdot)$ denotes the probability evaluated under the distribution $F$. In the following, we will mainly discuss properties related to $\alpha_+(F)$ and all the results can be generalized similarly to $\alpha_-(F)$. Both $\alpha_+(F)$ and $\alpha_-(F)$ take values in $[0,+\infty]$. For the right tail, $\alpha_+(F)=+\infty$ represents a thin tailed cdf such as the exponential distribution, and $\alpha_+(F)=0$ typically represents a super heavy tailed cdf such as the log-Pareto distribution (\citealt{CorRei09}).

The $\liminf$ used in the definition \eqref{tidef} is to pick up the heaviest part in the tail of $F$. The slowest possible decaying rate in the right tail of $F$ is roughly of order $O\left(x^{-\alpha_+(F)}\right)$ as $x\to +\infty$. Furthermore, if $F$ belongs to FMDA (i.e. it satisfies \eqref{righttail}), then the limit of the ratio $-\log \overline F(x) /\log x$ exists as $x\to +\infty$, and one can replace the $\liminf$ in \eqref{tidef} by $\lim$.

The definition of \eqref{tidef} can be viewed as a generalization of the usual tail index for distributions in FMDA. Frequentist estimators such as the Hill's estimator (\citealt{Hil75, DeHRes80}), the Pickands' estimator (\citealt{Pic75}) and their variations, are known to be asymptotically consistent for the tail index defined in \eqref{tidef} for certain restricted classes of distributions, such as FMDA. In general, it is unknown whether $\alpha_+(F)$ and $\alpha_-(F)$ defined by \eqref{tidef} can be consistently estimated from the data. However, this generalized notion of tail index in \eqref{tidef} is useful in describing the tail behavior of potentially complicated distributions drawn from Bayesian nonparametric priors.


\subsection{Bayesian Estimation and Posterior Consistency}\label{pcdef}


Let $\mathcal{F}$ be the set of all distributions that are absolutely continuous with respect to the Lebesgue measure on $\mathbb{R}$. Let $\mathscr{F}$ be the set of all density functions with respect to the Lebesgue measure on $\mathbb{R}$. Suppose that we observe a sample of i.i.d. data ${\bf X}^n = \{X_1,\ldots, X_n\}$ from the true distribution $F_0$ with the density $f_0$ on $\mathbb{R}$. Then in the Bayesian paradigm, we can impose a prior distribution on the probability density function $f \in \mathscr{F}$. For generality, here we will denote such a prior as $\Pi_n(\ud f)$, which explicitly allows the prior to depend on the sample size $n$. Equivalently, $\Pi_n$ is also a prior distribution over the set $\mathcal{F}$. Then the posterior distribution of $\Pi_n\left(\cdot|{\bf X}^n\right)$ evaluated at some measurable set $A\subseteq \mathscr{F}$ is
\begin{equation}\label{posterior}
\Pi_n\left(A|{\bf X}^n \right) = \frac{\int_A \prod_{i=1}^n f(X_i)\Pi_n(\ud f)}{\int_{\mathscr{F}} \prod_{i=1}^n f(X_i)\Pi_n(\ud f)}.
\end{equation}

To study properties related to the tail index, we define the following notion of tail index neighborhood.
\begin{definition}\label{tineibor}
For any distribution $F$ and $\epsilon>0$, the {\it $\epsilon-$(right) tail index neighborhood} of $F$ with $\alpha_+(F)\in [0,+\infty)$ is
$$B_{\alpha+}(F,\epsilon)\equiv \left \{H \in \mathcal{F}:~|\alpha_+(H)-\alpha_+(F)|<\epsilon\right\},$$
where $\alpha_+(\cdot)$ is defined in \eqref{tidef}. If $\alpha_+(F)=+\infty$, then the $\epsilon-$(right) tail index neighborhood of $F$ is defined as
$$B_{\alpha+}(F,\epsilon)\equiv \left \{H \in \mathcal{F}:~\alpha_+(H)=+\infty\right\}.$$
\end{definition}
The difference between the tail indices of two distributions used in Definition \ref{tineibor} is only a pseudometric, since different distributions can have the same tail index. In general, the topology induced by this pseudometric can be different from the weak topology generated by the weak convergence of probability measures. However, in the next proposition, we show that $B_{\alpha+}(F,\epsilon)$ is a Borel set on the space of all absolutely continuous distributions with respect to the Lebesgue measure on $\mathbb{R}$ associated to the weak topology. The proof is given in Appendix A.
\begin{prop}\label{borelset}
$B_{\alpha+}(F,\epsilon)$ is a Borel set on $\mathcal{F}$ under the weak topology, for any distribution $F$ and any $\epsilon>0$.
\end{prop}


Because the true tail index $\alpha_{0+} = \alpha_+(F_0)$ is unknown {\it a priori}, we hope that a distribution $F$ drawn from the posterior $\Pi_n\left(\cdot|{\bf X}^n \right)$ in \eqref{posterior} has a tail index $\alpha_+(F)$ sufficiently close to the truth $\alpha_{0+}$, as the sample size $n$ increases to infinity. This notion of asymptotics is usually stated as {\it consistency}.
\begin{definition}\label{pcti}
The posterior distribution $\Pi_n(\cdot|{\bf X}^n)$ is consistent for the (right) tail index if for any $\epsilon>0$, as $n\to \infty$,
$$\Pi_n \left(B_{\alpha+}^c(F_0,\epsilon) \big|{\bf X}^n\right) \to 0, \text{ in } P_{F_0}^{(n)} \text{ probability}.$$
\end{definition}
Definition \ref{pcti} is similar to the usual definition of posterior consistency for density estimation, but uses the tail index neighborhood in Definition \ref{tineibor}. It requires that the posterior probability assigns almost zero mass to distributions outside $\epsilon-$balls of $F_0$ as the sample size goes to infinity. On the other hand, although the weak consistency of density estimation is already well known for kernel mixture models \eqref{model2} below (see for example \citealt{Ghoetal99, Tok06, WuGho08}), posterior consistency of tail index does not follow directly from these results and requires further study, due to the non-equivalence between their topologies and neighborhoods.

\section{Tail Index of Location-Scale Mixture Models}\label{lsm}
In this section, we focus on a special case of Model \eqref{model1}, the location-scale mixture model
\begin{align}\label{model2}
f(x) = \int \frac{1}{\sigma}k\left(\frac{x-\mu}{\sigma}\right) \ud G(\mu,\sigma), \qquad G \sim \pi(G;\bxi),
\end{align}
where $k(\cdot)$ is a kernel density function and the parameter $\bt=(\mu,\sigma)$ consists of the location parameter $\mu$ and the scale parameter $\sigma$. We assume that the kernel $k(\cdot)$ has full support on $\mathbb{R}$.
Frequentist asymptotic properties of this model have been extensively studied in the Bayesian nonparametrics literature. Both weak and strong posterior consistency of Model \eqref{model2} have been discussed in \cite{Ghoetal99, Tok06, WuGho08}, etc. Theorem 3.3 of \cite{Tok06} established weak consistency of Model \eqref{model2} when the true density $f_0$ has a very thick polynomially decaying tail, with the tail index in $(0,1)$. However, in the following, we will show that weak consistency, and even strong consistency based on $L_1$ or Hellinger distance, is insufficient for meaningful Bayesian inference of the tail index. Surprisingly, for many commonly used priors $\pi(G;\bxi)$, the tail index of $F$ generated from Model \eqref{model2} can only take \emph{one single value}, implying that there is no possibility of identifying the correct tail index unless we know the true $\alpha_{0+}$ {\it a priori}.

For the MFM model (\citealt{RicGre97, GreRic01}), $f(x)$ in Model \eqref{model2} is specified as a finite mixture of $N$ components ($N\in \mathbb{Z}^+$), and a further prior distribution is imposed on $N$. In more details, the model is given as,
\begin{align}
\label{model3}
f\left(x \right) & = \sum_{i=1}^N \frac{w_i}{\sigma_i}k\left(\frac{x-\mu_i}{\sigma_i}\right), \nonumber \\
(\mu_i,\sigma_i)_{i=1}^N \Big|N &\stackrel{\text{iid}}{\sim} G_0(\mu,\sigma)\nonumber \\
(w_1,\ldots,w_N) \Big| N &\sim \text{Dirichlet}(a,\ldots,a), \;\text{for some}\; a>0\nonumber \\
N &\sim \pi(N) \text{ for } N =1,2,\ldots ,
\end{align}
The following theorem characterizes the tail index of a distribution $F$ generated by Model \eqref{model3}.

\begin{theorem}\label{mfmcounter}
Suppose that $G_0$ is a continuous distribution for $(\mu,\sigma)$. Then for any distribution $F$ with density $f$ drawn from Model \eqref{model3}, the range of $\alpha_+(F)$ is almost surely a singleton. In other words, almost surely all $F$'s drawn from the MFM model have the same tail index.
\end{theorem}
In the finite mixture model given in \eqref{model3}, the tail indices of different $F$'s are all the same, since all of them are finite mixtures and their tail indices are solely determined by the tail heaviness of the kernel $k(\cdot)$. A heavy tailed kernel will only make the tails of $F$ heavy, but not be able to generate varying tail heaviness. This limitation immediately indicates that we cannot obtain any meaningful posterior consistency in terms of tail index.

We now investigate the more complicated example where $G(\mu,\sigma)$ has a nonparametric NRMI prior. In the theorems to follow, we adopt similar NRMI notations as those in \cite{Lijetal07}, \cite{Jametal09}, and \cite{Baretal13}. We consider a completely random measure $\widetilde H$, with $\widetilde H(x) = \sum_{i\geq 1} \widetilde J_i\delta_{X_i}(x)$ for $x\in \mathbb{R}$ such that $\{X_i\}_{i\geq 1}$ and nonnegative $\{\widetilde J_i\}_{i\geq 1}$ are independent sequences of random variables, ignoring jumps at nonrandom positions. The joint distribution of $\{\widetilde J_i\}_{i\geq 1}$ and $\{X_i\}_{i\geq 1}$ is characterized by the L\'evy intensity $\nu(\ud v,\ud x)$ through the Laplace transformation of $\widetilde H$ (for $s>0$):
\begin{align*}
E\left[e^{-s\widetilde H(A)}\right] = \exp\left\{-\int_{\mathbb{R}^+\times A} \left(1-e^{-sv}\right)\nu (\ud v,\ud x)\right\}, \quad \text{for any } A\subseteq \mathbb{R}.
\end{align*}
We consider the homogenous NRMI where the L\'evy intensity can be factorized as $\nu (\ud v,\ud x) = \rho(\ud v) H_0(\ud x)$. $\rho(\ud v)$ is the L\'evy intensity for the nonnegative masses $\{\widetilde J_i\}_{i\geq 1}$, and $\{X_i\}_{i\geq 1}$ are independent draws from the nonatomic probability measure $H_0$, also called the ``base measure". Then a NRMI $H$ is defined as $H(x) =\sum_{i\geq 1} J_i\delta_{X_i}(x)$ with $J_i=\widetilde J_i/ \sum_{i\geq 1}\widetilde J_i$ for any $x\in \mathbb{R}$. For all the theorems in this section, we assume that $\rho(\ud v)$ satisfies $\int_0^{\infty}\rho(\ud v)=+\infty$ and $\int_0^{\infty}(1-e^{-v})\rho(\ud v)<+\infty$ which guarantees that $0<\sum_{i\geq 1} \widetilde J_i<+\infty$ almost surely and the NRMI $H$ is well defined; see equation (2.3) of \citet{FavTeh13}.

The following theorem will be used as a fundamental tool in studying the tail behavior of a NRMI.
\begin{theorem}\label{nrmcounter1}
Suppose $H$ is a homogeneous NRMI with the L\'evy intensity measure $\rho(\ud v)H_0(\ud x)$ for $v\in \mathbb{R}^+,x\in \mathbb{R}$ where $H_0$ is a continuous probability measure on $\mathbb{R}$. Let $\Psi(s) = \int_0^{+\infty}(1-e^{-sv})\rho(\ud v)$ and let $\Psi^{-1}$ be the inverse function of $\Psi$. Then \\
\noindent (i) If there exists a function $h_{\gamma}$ defined as
\begin{align}\label{hgamma}
&h_{\gamma}(x) = \frac{\log |\log x|}{\Psi^{-1}\left(\frac{\gamma \log |\log x|}{x}\right)},
\end{align}
with $\gamma>1$ for $x\in (0,1/e)$, such that\\
$\liminf_{x\to +\infty} -\log h_{\gamma}(\overline H_0(x))/\log x = 0$, then $\alpha_+(H)=0$ a.s.\\
\noindent (ii) If there exists a function $h$ such that\\
\hspace*{3mm} (a) $h(x)$ is locally convex in $x\in [0,\epsilon)$ for some small $\epsilon>0$;\\
\hspace*{3mm} (b) $\int_0^\epsilon \rho[h(x),+\infty)\ud x < +\infty$;\\
\hspace*{3mm} (c) $\liminf_{x\to +\infty} -\log h(\overline H_0(x))/\log x = +\infty$;\\
then $\alpha_+(H)=+\infty$ a.s.
\end{theorem}
The theorem follows from \cite{Fri67} and \cite{FriPru71}; see Proposition \ref{Friprop} and the subsequent proof of Theorem \ref{nrmcounter1} in Appendix A. Proper choices of the functions $h_{\gamma}$ in (i) and $h$ in (ii) will lead to sharp lower and upper bounds for the tail index of a NRMI $H$.


The next theorem describes how these bounds for a NRMI can be used to characterize the tail behavior of a mixture density drawn from Model \eqref{model2}.


\begin{theorem}\label{nrmcounter2}
Suppose in Model \eqref{model2}, $G(\cdot,\cdot)$ is a homogeneous NRMI with L\'evy intensity measure $\rho(\ud v)G_0(\ud \mu,\ud\sigma)$ for $v\in \mathbb{R}^+$ and $(\mu,\sigma) \in \mathbb{R}\times \mathbb{R}^+$, where $G_0(\mu,\sigma)$ is a continuous cdf on $\mathbb{R}\times \mathbb{R}^+$. Let $G_{0,\mu}$ and $G_{0,\sigma}$ be the marginal distributions of $G_0(\mu,\sigma)$ for $\mu$ and $\sigma$, respectively. Assume that $G_{0,\mu}$ is symmetric about zero. If both $G_{0,\mu}$ and $G_{0,\sigma}$ satisfy either (i) or (ii) in Theorem \ref{nrmcounter1}, i.e. we can replace $H_0$ in either (i) or (ii) of Theorem \ref{nrmcounter1} by $G_{0,\mu}$ and $G_{0,\sigma}$, then for any distribution $F$ with density $f$ drawn from Model \eqref{model2}, the range of $\alpha_+(F)$ as defined in \eqref{tidef} is almost surely a singleton.
\end{theorem}

Theorem \ref{nrmcounter2} indicates that the tail indices of all distributions $F$ drawn from Model (\ref{model2}) are almost surely the same, if each of the two marginals $G_{0,\mu}$ and $G_{0,\sigma}$ satisfies either (i) or (ii) in Theorem \ref{nrmcounter1}. Again this indicates that there is no meaningful posterior consistency for tail index, by similar arguments after Theorem \ref{mfmcounter}. Theorem \ref{nrmcounter2} and its proof also lead to two other interesting implications. First, if the conditions for the two marginals of the base measure hold, then the tail index of $F$ only depends on the tail indices of the two marginals, but not on the full joint distribution $G_0(\mu,\sigma)$. Second, whether $\alpha_+(F)$ is the same for all $F\sim \Pi_n(\cdot|{\bf X}^n)$ does not depend on the tail behavior of the kernel $k(\cdot)$, even if $k(\cdot)$ is a heavy tailed kernel.

\begin{remark}
The assumption of symmetric $G_{0,\mu}$ is only used as a sufficient condition for the case where $G_{0,\mu}$ satisfies (ii) of Theorem \ref{nrmcounter1} and $G_{0,\sigma}$ satisfies (i) of Theorem \ref{nrmcounter1}, in other words, the case where $G_{\mu}$ has thin left and right tails and $G_{\sigma}$ has a super heavy right tail. The assumption of symmetric $G_{0,\mu}$ is not necessary for the conclusions of Theorem \ref{nrmcounter2} to hold when both $G_{\mu}$ and $G_{\sigma}$ are thin tailed, and when $G_{\mu}$ has a super heavy right tail. Details of the proof can be found in Appendix A.
\end{remark}

\begin{remark}
The proof of Theorem \ref{nrmcounter2} also relies on the moment techniques in Lemma A.1--A.3 in Appendix A, which relate the tail index of $F$ to the moments of $F$, and subsequently the moments of the kernel $k(\cdot)$ and the mixing distribution $G$. As a side product of this proof, we recovered the famous Breiman lemma in \cite{Bre65} about scale mixtures with heavy tailed mixing measures. Suppose in Model \eqref{model2} we only have the scale mixture $f(x)=\int \sigma^{-1}k(x/\sigma) \ud G(\sigma)$ and $G$ has tail index $\alpha_+(G)\in (0,+\infty)$. The Breiman lemma says that if the kernel $k(\cdot)$ has a tail index larger than $\alpha_+(G)$, i.e. it has a thinner tail than $G$, then the mixture $f(x)$ is also heavy tailed with tail index $\alpha_+(G)$. This is an immediate result of our Lemma A.3.
\end{remark}

We make the tail conditions on $G_{0,\mu}$ and $G_{0,\sigma}$ more concrete for the special cases of Dirichlet process (DP) and normalized generalized Gamma process (NGGP, \citealt{Lijetal07, Jametal09, LijPru10, Baretal13}) mixture models.  It turns out that there is a large class of measures that satisfy the condition (i) or (ii) in Theorem \ref{nrmcounter2}, including both thin tailed distributions and heavy tailed distributions.

\begin{theorem}\label{dpcounter}
Suppose in Model \eqref{model2}, $G(\cdot,\cdot)\sim \DP (a,G_0(\mu,\sigma))$ with $a>0$ and $G_0(\mu,\sigma)$ a continuous cdf on $\mathbb{R}\times \mathbb{R}^+$. Assume that $G_{0,\mu}$ is symmetric about zero. Consider the following two conditions for a generic distribution $H_0$ on $\mathbb{R}$:\\
\noindent (i) $\limsup_{x\to +\infty}\overline H_0(x)\cdot \left(\log x /\log\log x\right) = +\infty$,\\
\noindent (ii) $\limsup_{x\to +\infty}\overline H_0(x) \cdot \left[(\log x)\cdot (\log\log x)^{\delta}\right] = 0$ for some $\delta>1$.\\
If both $G_{0,\mu}$ and $G_{0,\sigma}$ satisfy either one of the conditions (i) and (ii), i.e. we can replace $H_0$ in either (i) or (ii) by $G_{0,\mu}$ and $G_{0,\sigma}$, then for any distribution $F$ with density $f$ drawn from Model \eqref{model2}, the range of $\alpha_+(F)$ as defined in \eqref{tidef} is almost surely a singleton.
\end{theorem}
The proof of Theorem \ref{dpcounter} involves the tail behavior of a DP, which has been studied in \cite{DosSel82}. Conditions (i) and (ii) correspond to conditions (i) and (ii) in Theorem \ref{nrmcounter2}. As a result of the theorem, most distributions $G_{0,\mu}$ and $G_{0,\sigma}$ with either heavier or thinner tails than $1/\log x$ will lead to a single value of tail index for all $F$'s in the DP mixture model, and therefore the posterior cannot estimate the truth $\alpha_{0+}$ consistently. For example, in the popular DP mixture of normals (\citealt{EscWes95}), the marginal distributions of the base measure for $\mu$ and $\sigma^2$ are the Student's $t$ distribution and the inverse gamma distribution, both of which have much thinner tails than $1/\log x$. Therefore, the Bayesian posterior with the normal-inverse gamma prior for DP mixture of normals cannot consistently estimate the tail index. In contrast, Theorem 3.3 of \cite{Tok06} has shown that such a normal-inverse gamma base measure is sufficient for posterior weak consistency, even if the true density is heavy tailed with a tail index in $(0,1)$. This implies that the conditions required for consistent estimation of the tail index are more stringent than those for usual weak and strong posterior consistency. We emphasize again that the kernel here plays an inconsequential role due to Theorem \ref{nrmcounter2}, regardless of its tail thickness.

An important implication of Theorem \ref{dpcounter} is that the bounds in (i) and (ii) are not far from each other. As a result, not many distributions have been left out by (i) and (ii). Basically, only those base measures that decay at a similar rate to $1/\log x$ are not covered by the conditions (i) and (ii). As a result, the only combination that is not covered by Theorem \ref{dpcounter} is the case where both $G_{0,\mu}$ and $G_{0,\sigma}$ decay at rates similar to $1/\log x$. When this happens, the tail index of $F$ drawn from Model \eqref{model2} can possibly vary in $[0,+\infty]$. In this case, whether the posterior consistency of tail index holds or not remains unknown.


The next theorem shows a similar posterior behavior for the general $\NGGP$ mixture model, denoted by $\NGGP(a,\kappa,\tau,G_0(\mu,\sigma))$. Its L\'evy intensity measure is given by $\rho(\ud v)dG_0(\mu,\sigma) = \frac{a}{\Gamma(1-\kappa)}v^{-\kappa-1}e^{-\tau v} \ud v dG_0(\mu,\sigma)$, where $a>0$, $\kappa\in [0,1)$ and $\tau>0$. The NGGP class includes most of the discrete random probability measures in the Bayesian nonparametric literature. For example, the class includes DP as $\NGGP(a,0,1,G_0)$, the normalized-inverse Gaussian process as $\NGGP(1,1/2,\tau,G_0)$, and the N-stable process as $\NGGP(1,\kappa,0,G_0)$ as special cases. See \cite{Lijetal07} and \cite{Baretal13} for discussions. The cases of $\kappa=0$ (DP) and $\kappa>0$ are different in nature, so the conclusion of Theorem \ref{nggpcounter} is also different from Theorem \ref{dpcounter}.

\begin{theorem}\label{nggpcounter}
Suppose in Model \eqref{model2}, $G(\cdot,\cdot)\sim \NGGP(a,\kappa,\tau,G_0(\mu,\sigma))$ with $a>0$, $\kappa\in (0,1)$, $\tau\geq 0$ and $G_0(\mu,\sigma)$ is a continuous cdf on $\mathbb{R}\times \mathbb{R}^+$. Assume that $G_{0,\mu}$ is symmetric about zero. Consider the following two conditions for a generic distribution $H_0$ on $\mathbb{R}$:\\
\noindent (i) $\limsup_{x\to +\infty}\overline H_0(x) \cdot x^{\delta} = +\infty$ for all $\delta>0$,\\
\noindent (ii) $\limsup_{x\to +\infty}\overline H_0(x) \cdot  x^{\delta} =0 $ for all $\delta>0$.\\
If both $G_{0,\mu}$ and $G_{0,\sigma}$ satisfy either one of the conditions (i) and (ii), i.e. we can replace $H_0$ in either (i) or (ii) by $G_{0,\mu}$ and $G_{0,\sigma}$, then for any distribution $F$ with density $f$ sampled from Model \eqref{model2}, the range of $\alpha_+(F)$ as defined in \eqref{tidef} is almost surely a singleton.
\end{theorem}
Similar to Theorem \ref{dpcounter}, here we also provide two conditions for the tail decaying rates of $G_{0,\mu}$ and $G_{0,\sigma}$, where (i) gives heavier than polynomial tails and (ii) gives thinner than polynomial tails. The gap between the base measures that satisfy (i) or (ii) in the current theorem is now larger than that in the DP case, but the theorem still has ruled out many possibilities for consistent estimation of tail index. For example, when both $G_{0,\mu}$ and $G_{0,\sigma}$ have exponentially decaying tails, the tail index generated from the posterior of a NGGP is always the same as the tail index of the kernel $k(\cdot)$ (see the proof of Theorem \ref{nrmcounter2} in Appendix A). It remains unknown how the tail indices of $F$ from a NGGP mixture model behave in the posterior when at least one of $G_{0,\mu}$ and $G_{0,\sigma}$ have a polynomially decaying tail.

\section{Sufficient Conditions for Tail Index Consistency}
\label{suff}
\subsection{Schwartz's Theorem for Posterior Consistency}
In this section, we provide a series of conditions that guarantee the posterior consistency of tail index for the most general model $f\sim \Pi_n$. These conditions are built on the classic Schwartz's argument in \cite{Sch65} for posterior consistency, and therefore they are simple and intuitive. We will then demonstrate the application of these sufficient conditions on Model \eqref{model1} using the Pareto kernel in Section \ref{example}.

The definition of tail index in \eqref{tidef} applies to any distribution but may be too general so that no consistent frequentist estimator exists. Therefore, we will limit our scope to those priors that only generate candidate distributions from the class of FMDA, i.e. distributions that satisfy \eqref{righttail}. These distributions have a well defined tail index, i.e. we can replace all the $\liminf$ in \eqref{tidef} by $\lim$. Throughout the entire Section \ref{suff}, we assume that the true distribution has a tail index $\alpha_{0+}\in (0,+\infty)$, and the prior $\Pi_n$ satisfies Condition (PT).
\vspace{.4cm}

\noindent (PT) For almost surely all $F\sim \Pi_n$, $F$ satisfies the relation \eqref{righttail} with $\alpha_+(F)\in (0,+\infty)$ and a slowly varying function $L_F$, and its right tail index is given by \eqref{tidef} with all $\liminf$ replaced by $\lim$.
\vspace{.4cm}

The Schwartz consistency theorem relies on two key conditions: the Kullback-Leibler (KL) support of the prior, and the existence of a uniformly consistent test. For two distributions $F_1$ and $F_2$ (with densities $f_1$ and $f_2$), let the KL divergence between $F_1$ and $F_2$ be $KL(F_1,F_2)\equiv E_{F_1}\log (f_1/f_2)$. Define the $\epsilon-$KL neighborhood of the true distribution $F_0$ as $\mathcal{K}(F_0,\epsilon)\equiv \{F \in \mathcal{F}: KL(F_0,F) <\epsilon\}$. The condition on the KL support of the prior is stated as follows:
\vspace{.4cm}

\noindent (KL) The true distribution $F_0$ is in the KL support of $\Pi_n$, if for any $\epsilon>0$, \\
$\liminf_{n\to \infty} \Pi_n (\mathcal{K}(F_0,\epsilon))>0$.
\vspace{.4cm}

We allow the prior $\Pi_n$ to depend on the sample size $n$, since this can be conveniently incorporated into the standard posterior consistency argument (see Section 5 of \citealt{Ghoetal99}). It is well known that the condition (KL) implies weak consistency, and is therefore a very basic requirement for useful Bayesian models.

The other condition required in the Schwartz consistency theorem is the existence of uniformly consistent tests. For our purpose, we need a test for tail index that is able to separate $F_0$ from all the distributions \emph{outside a tail index neighborhood of $F_0$}. A set $\mathcal{F}_n$ with large prior probability (called ``sieve") helps when $\Pi_n$ has a non-compact support and the uniform test can be found on a sufficiently large set.
\vspace{.4cm}

\noindent (UT) Uniform testing condition: There exists a test $\Phi_n\equiv \Phi_n(X_1,\ldots,X_n)$ and a sieve $\mathcal{F}_n$ such that \\
\noindent (i) $\Pi_n(\mathcal{F}_n^c)\leq e^{- bn}$ for some constant $b>0$; \\
\noindent (ii) For any $\epsilon>0$, as $n\to\infty$,
\begin{align}\label{test3}
&E_{F_0}\Phi_n \to 0, \qquad \sup_{F\in B^c_{\alpha+}(F_0,\epsilon)\cap \mathcal{F}_n} E_{F}(1-\Phi_n)\to 0.
\end{align}

Based on Schwartz's consistency theorem, one can show posterior consistency of tail index under the conditions (KL) and (UT).
\begin{theorem}\label{schthm}
If both (KL) and (UT) hold true, then the posterior distribution $\Pi_n(\cdot|{\bf X}^n )$ is consistent for the (right) tail index.
\end{theorem}
The proof follows the same thread as the usual weak consistency (see for example \citealt{Ghoetal99}, \citealt{GhoRam03}) and is therefore omitted. Note that the uniform test in (UT) can be made exponentially fast by an argument using the Hoeffding's inequality (Theorem 2 of \citealt{Ghoetal99}, Proposition 4.4.1 of \citealt{GhoRam03}). However, a key unanswered question is whether such a uniformly consistent test $\Phi_n$ for tail index exists. One cannot directly apply the Le Cam theory because $\Phi_n$ will depend on the new tail index neighborhood of $B_{\alpha+}(F_0,\epsilon)$ and the pseudometric about tail index difference. We instead proceed in a constructive way and pursue sufficient conditions for (UT) to hold.

\subsection{Existence of Tests}

In the representation \eqref{righttail} for a generic distribution $F\sim \Pi_n$, let $h_F(x) = xL'_F(x)/L_F(x)$ and hence $L_F(x) = L_F(x_0)\exp\left(\int_{x_0}^x \frac{h_F(t)}{t} \ud t\right)$ for some fixed $x_0$. Alternatively, $h_F(x)$ can be written as
$$h_F(x) = \alpha_+(F) - \frac{xf(x)}{\overline F(x)}.$$
For any given $F$ from FMDA, the von-Mises theorem (see e.g. Proposition 2.1 of \citealt{Beietal04}) says that
$$\lim_{x\to +\infty} \frac{xf(x)}{\overline F(x)} = \alpha_+(F),$$
i.e. $\lim_{x\to +\infty}h_F(x)=0$. Bounding the magnitude of $h_F(x)$ is crucial in showing the existence of uniform tests for $\alpha_+(F)$. In the Bayesian framework, $h_F(x)$ with $F\sim \Pi_n$ needs to be controlled in a uniform way on a sieve with large prior probability. In light of this, we have the following theorem on the existence of tests. Throughout the rest of the paper, for two positive sequences $\{x_n\}$ and $\{y_n\}$ that depend on the sample size $n$, $x_n \prec y_n$ means $x_n=o(y_n)$, $x_n\succ y_n$ means $y_n = o(x_n)$, $x_n \preceq y_n$ means $x_n = O(y_n)$, and $x_n\succeq y_n$ means $y_n=O(x_n)$.

\begin{theorem}\label{utest1}
Suppose that $\alpha_{0+}\in (0,+\infty)$, and (PT) holds. In addition, suppose the following conditions hold:\\
\noindent (i) There exist finite constants $x_0 \geq e$ and $c_L\in (0,1)$, such that for all sufficiently large $n$, $L_F(x_0)\geq n^{-c_L}$ uniformly for all $F\in \mathcal{F}_{1n}$, where $\mathcal{F}_{1n}$ is a sieve satisfying $\Pi_n(\mathcal{F}_{1n}^c) < e^{-c_1n}$ for some constant $c_1>0$;\\
\noindent (ii) There exists an envelope function $\overline h_n(x)=B_n (\log x)^{-(1+\tau_n)}$ for some positive $n$-dependent sequences $B_n$ and $\tau_n$, such that for all sufficiently large $n$, $|h_F(x)|\leq \overline h_n(x)$ for all $F\in \mathcal{F}_{2n}$ and all $x\geq x_0$, where $\mathcal{F}_{2n}$ is a sieve satisfying $\Pi_n(\mathcal{F}_{2n}^c) < e^{-c_2n}$ for some constant $c_2>0$; \\
\noindent (iii) The prior $\Pi_n$ satisfies $\Pi_n(\mathcal{F}_{3n}^c) < e^{-c_3 n}$ for some constant $c_3>0$ for $\mathcal{F}_{3n}=\{F\in \mathcal{F}:\alpha_+(F) \leq \overline \alpha_n \}$ and some sequence $1\prec \overline \alpha_n \prec \log n$, for all sufficiently large $n$;\\
\noindent (iv) $B_n$, $\tau_n$ and $\overline \alpha_n$ satisfy $1\preceq B_n \prec \min(\overline \alpha_n ^{-1} \log n, \tau_n\log n)$ and $\tau_n \preceq 1$;\\
then (UT) holds.
\end{theorem}

The proof of Theorem \ref{utest1} uses a recently proposed tail index estimator in \cite{CarKim15} defined as
\begin{align}\label{alphatest}
&\hat \alpha_{s_n} = \log(\hat p_{s_n}) - \log(\hat p_{s_n+1}),
\end{align}
where $\hat p_{s_n} = n^{-1} \sum_{i=1}^n I(X_i>e^{s_n})$ and $s_n$ is taken as a positive sequence that satisfies $B_n \prec s_n \prec \overline \alpha_n ^{-1} \log n$ (see the proof of Theorem \ref{utest1} in Appendix A). Such a sequence $s_n$ exists given Condition (iv) in Theorem \ref{utest1}. \cite{CarKim15} has shown that when $\alpha_{0+}\in (0,+\infty)$, $\hat \alpha_{s_n}$ is a consistent estimator of $\alpha_+(F)$ for $F$ from various classes of distributions, such as the first order and the second order approximately Pareto distributions. \cite{CarKim15} has also given the explicit choice of $s_n$ (as well as a data-dependent version) such that $\hat \alpha_{s_n}$ converges at a minimax rate to $\alpha_+(F)$ for a certain class of distributions (adaptively). Therefore, a test for $H_0:\alpha_+(F)=\alpha_{0+}$ can be $\Phi_n=I\left(|\hat \alpha_{s_n}-\alpha_{0+}|>\epsilon\right)$ given some $\epsilon>0$. For our purposes, it is easier to work with $\hat \alpha_{s_n}$ than the Hill's estimator.

Conditions (i)-(iv) are sufficient for the existence of such tests. Among them, (i) and (ii) are mainly intended to control the slowly varying function $L_F$, where we allow exceptions on sets with exponentially small prior probabilities. The choice of $x_0\geq e$ is mainly for convenience since $\log x>1$ for all $x\geq x_0$. Alternatively, one can replace it with any finite $x_0\in \mathbb{R}$ and modify the definition of logarithm function with a shift accordingly. In (ii) we specify the envelope function $\overline h_n(x)$ to be decaying in the logarithm of $x$. In the frequentist tail index literature, such control over the exponent in a slowly varying function has appeared in \cite{Dre98} and \cite{Dre01} for showing minimax rates in certain classes of distributions. The logarithmically decaying $\overline h_n(x)$ is not restrictive because we allow $B_n\to \infty$ and $\tau_n\to 0$ as $n\to \infty$. As an envelop function, it also includes all $h_F(x)$ that decays polynomially in $x$.




Condition (iii) restricts the largest possible tail index on a large sieve, but the sieve will eventually cover the true $F_0$ as the sample size $n$ increases. Condition (iv) determines the choice of $B_n,\tau_n$ in (ii) and $\overline\alpha_n$ in (iii). For posterior consistency, we only require the existence of such sequences $B_n,\tau_n,\overline\alpha_n$. Conditions (i)-(iv) will be verified for Pareto mixtures in Section \ref{example}.

\begin{remark}
We would like to emphasize that in our Bayesian setup, the class of distributions for which $\hat \alpha_{s_n}$ in \eqref{alphatest} gives a uniform test depends on the conditions on the prior, for example Conditions (i)-(iv) in Theorem \ref{utest1}. These conditions impose restrictions on the class of distributions and densities that can be consistently fitted by our posterior \eqref{posterior} in the sense of weak consistency. In fact, they can result in a relatively smaller KL support of the prior, which may only include a subclass of FDMA. This is partly due to the basic requirement that our Baysian posterior should achieve consistency for both fitting the density and fitting the tail index at the same time. Although in general it is difficult to describe exactly which distributions are included in the prior KL support given those conditions in Theorem \ref{utest1}, we will shed light on this for the example of Pareto mixtures in Theorem \ref{cmkl} in Section \ref{example}.
\end{remark}

The following theorem is a consequence of Theorem \ref{schthm} and Theorem \ref{utest1}.
\begin{theorem}\label{post1}
(Posterior Consistency of Tail Index) Under all assumptions of Theorem \ref{utest1} and (KL), the posterior distribution $\Pi_n (\cdot|{\bf X}^n )$ is consistent for the (right) tail index.
\end{theorem}

%
%

\subsection{Example of Consistency: Mixtures of Paretos}\label{example}
The failure of tail index consistency in Section \ref{lsm} is partly due to the structure of the location-scale mixture model \eqref{model2}, in which we have no control over how the mixing measure $G(\mu,\sigma)$ affects the tail index of the mixture distribution. A possible remedy is to introduce an explicit mixture on the tail index parameter. An example of this type is the DPM of Paretos used in \cite{Tre08}. In this section, we study the mixture of simple Pareto distributions with kernel density $k(x;\alpha)=\alpha x^{-(\alpha+1)}$ whose support is $[1,+\infty)$. We will take the mixing measure from a homogenous NRMI prior, such as DP and NGGP. Because a general discrete mixture distribution takes the form $\overline F(x)=\sum_{i=1}^{\infty} w_i x^{-\alpha_i}$, the right tail index is $\alpha_+(F)=\inf\{\alpha_1,\alpha_2,\ldots\}$. To make this tail index more explicit, in the following Bayesian model, we are going to first pick $\alpha_1$ as the tail index of $F$ together with its weight $w_1$, and then draw the other $\alpha_i$ and their weights $w_i$ ($i=2,3,\ldots$) from a mixture model conditional on $\alpha_1$ and $w_1$. In this way we can guarantee that $\alpha_i>\alpha_1$ for all $i\geq 2$ such that we can conveniently control the behavior of $\alpha_+(F)$ through $\alpha_1$. The model is specified as follows.
\begin{align}\label{parmix}
f(x) \Big| \alpha_1, w_1, H & = w_1 k\left(x;\alpha_1\right) + (1-w_1)\int k(x;\alpha) \ud H(\alpha), \nonumber \\
\alpha_1 & \sim G_{\alpha}\cdot I_{\left[0,\overline \alpha_n\right]}, \text{ supp}(G_{\alpha}) = \left[0,+\infty\right),~ G_{\alpha}\text{ has no point mass at zero}, \nonumber \\
w_1 & \sim G_w \cdot I_{\left[\underline w_n,1\right]}, \text{ supp}(G_w) = \left[0,1\right], \nonumber \\
H_1 & \sim \Pi(H_1;\bxi,H_0), \text{ supp}(H_0)=[0,+\infty), ~ H_0 \text{ has no point mass at zero},\nonumber \\
H(\alpha) & = H_1(\alpha-\alpha_1), \text{ for any } \alpha>\alpha_1.
\end{align}
The notation ``supp" stands for the support of a distribution. For a generic distribution $G$ and a set $A$, $G\cdot I_A$ denotes the renormalized probability distribution of $G$ truncated to the set $A$. The density $f$ has two mixing components. The first component $w_1 k\left(x;\alpha_1\right)$ explicitly controls the tail index of $F$, and the second component is a general mixture of Paretos. $\alpha_1$ in the first component determines $\alpha_+(F)$, and is drawn from $G_{\alpha}$ truncated to $[0,\overline \alpha_n]$. Here the deterministic positive sequences $\underline w_n$ and $\overline \alpha_n$ satisfy that $\underline w_n\to 0$ and $\overline \alpha_n\to +\infty$ as $n\to \infty$, so asymptotically the supports of $w_1$ and $\alpha_1$ covers any number in $(0,1]$ and $\mathbb{R}^+$. The second component in the mixture involves a mixing probability measure $H$, which is drawn from a prior $\Pi$. $\bxi$ contains all the hyperparameters of $\Pi$, such as the parameter $a$ in a DP and the parameters $a,\kappa,\tau$ in a NGGP. Given the value of $\alpha_1$, $H$ is a right-shifted version of the distribution $H_1$ drawn from the prior $\Pi$. For the ease of presentation, we assume that $G_{\alpha},G_w$ and $\Pi$ do not depend on $n$.






The deterministic sequences $\underline w_n$ and $\overline \alpha_n$ introduced here are mainly designed to separate the leading component $w_1k(x;\alpha_1)$ from the other mixing components, such that the sufficient conditions in Theorem \ref{utest1} are satisfied. In particular, condition (PT) can be conveniently verified for Model \eqref{parmix} with the help from the leading component. $\overline \alpha_n$ is used such that $\alpha_1$ has an increasingly large support and meanwhile Condition (iii) of Theorem \ref{utest1} is satisfied. In fact, the way of isolating the leading Pareto component in Model \eqref{parmix} is similar to some well studied nonparametric classes of distributions in the frequentist tail index literature, such as the Hall and Welsh class (\citealt{HalWel85, CarKim15, BouTho15}) that satisfies $\left|\overline F(x)-Cx^{-\alpha}\right|\leq C'x^{-\alpha(1+\beta)}$ for $\alpha,\beta,C,C'>0$.



A function $g$ on the interval $I$ is called completely monotone if the $m$th derivative of $g$ satisfies $(-1)^mg^{(m)}(x)\geq 0$ for all $m\in \mathbb{Z}^+$. Let
\begin{align*}
\mathcal{CM}_{e} =& \left\{F: \text{supp}(F)=[1,+\infty),\overline F(e^t) \text{ is completely monotone on }t\in [0,+\infty)\right\},\\
\mathcal{P}_{2} =& \Big\{F: \text{supp}(F)=[1,+\infty), \overline F(x) = Cx^{-\alpha} + O(x^{-(1+\beta)\alpha}),\\
&\text{ for some constant } \alpha>0,\beta>0, C>0\Big\},
\end{align*}
where $\mathcal{P}_2$ is the class of second-order Pareto distributions. We can characterize the class of distributions described by Model \eqref{parmix}.
\begin{theorem}\label{cmkl}
Suppose in Model \eqref{parmix}, $\underline w_n\to 0$ and $\overline\alpha_n\to+\infty$ as $n\to\infty$. If $F\in \mathcal{CM}_e\cap \mathcal{P}_2$ and the prior $\Pi(H;\xi,H_0)$ is a homogeneous NRMI, then $F$ is in the KL support of Model \eqref{parmix}.
\end{theorem}

The KL support of Model \eqref{parmix} is related to the class of completely monotone functions. This is not surprising because the mixtures of Paretos are related to the mixtures of exponential distributions by the transformation $x=e^t$ in the Pareto kernel $k(x;\alpha)$. The KL support of the mixtures of exponentials includes the class of completely monotone functions (Theorem 16 in \citealt{WuGho08}), by the Hausdorff-Bernstein-Widder theorem. In fact, it is proved in Lemma S.1 in the supplementary material that any distribution $F$ from $\mathcal{CM}_e\cap \mathcal{P}_2$ has a density with a similar form to that in Model \eqref{parmix}.

The following theorem imposes further conditions on $\underline w_n, \overline \alpha_n$ and the prior $G_{\alpha},G_w,\Pi$, such that Model \eqref{parmix} achieves posterior consistency of tail index.

\begin{theorem}\label{parthm}
Suppose the following conditions hold for Model \eqref{parmix}:\\
(i) $F_0\in \mathcal{CM}_e \cap \mathcal{P}_{2}$;\\
(ii) The prior $\Pi$ on the mixing measure $H$ satisfies one of the following conditions:\\
\hspace*{3mm} (a) $\Pi$ is $\DP(a,H_0)$ where $a>0$ and $H_0$ is a probability distribution on $\mathbb{R}^+$, and there exist positive constants $0<c_1<1,D_1>0,d_1>0$, such that $H_0(x)\leq D_1[\log (1/x)]^{-(1+d_1)}$ for all $x\in (0,c_1)$; \\
\hspace*{3mm} (b) $\Pi$ is $\NGGP(a,\kappa,\tau,H_0)$ where $a>0$, $\kappa\in (0,1)$, $\tau>0$ and $H_0$ is a probability distribution on $\mathbb{R}^+$, and there exist positive constants $0<c_2<1,D_2>0,d_2>0$, such that $H_0(x)\leq D_2 x^{1+d_2}$ for all $x\in (0,c_2)$;\\
(iii) $1\prec \overline \alpha_n\prec \log n$, $\overline \alpha_n /\log n \prec \underline w_n \prec 1$;\\
then the posterior distribution $\Pi_n(\cdot|{\bf X}^n )$ of Model \eqref{parmix} is consistent for the tail index.
\end{theorem}
Condition (ii) in Theorem \ref{parthm} requires sufficient decay for the base measure $H_0$ near zero, though the decaying rate could be different for a DP prior and a NGGP prior. For a NGGP prior, the decaying rate of $H_0(x)$ near $x=0$ needs to be in polynomials of $x$, while the rate for a DP prior can be slower, in polynomials of $\log (1/x)$ for $x$ close to zero. This is due to the difference in the tail behavior of DP and NGGP. Condition (iii) describes the orders of $\underline w_n$ and $\overline\alpha_n$. They can be taken as, for example, $\underline w_n = (\log n)^{-1/3}$ and $\overline \alpha_n = (\log n)^{1/2}$.


\begin{remark} The densities in $\mathcal{CM}_e \cap \mathcal{P}_{2}$ always have nonnegative mixing coefficients, since $w_1>0$ and $H$ is a probability measure. As a result, the KL support of Model \eqref{parmix} includes mixtures such as $\overline F(x) = \frac{1}{2x} + \frac{1}{2x^2}$, but also has excluded some other mixtures of Paretos, such as $\overline F(x) = \frac{2}{x}-\frac{1}{x^2}$ in which some components may have negative coefficients. To enlarge the KL support of Model \eqref{parmix} and allow negative mixing coefficients, the mixing measure can be characterized as a bounded signed measure $w_1\delta_{\alpha_1}+(1-w_1)H = H_+-H_-$, where $\delta_a$ denotes the Dirac measure at $a$. Similar priors to those in Model \eqref{parmix} can be imposed on both $H_+$ and $H_-$ and they need further restrictions to guarantee that the density $f$ is nonnegative.  For example, if $\overline F(x) = \frac{2}{x}-\frac{1}{x^2}$, then $H_+=2\delta_1$ and $H_-=\delta_2$. According to Theorem 4.3 of \cite{Wat60}, the Pareto kernel mixture representation by using bounded signed mixing measure includes all distributions $F$ that satisfy $\sum_{i=1}^{\infty}\left|\overline F^{(i)}(e^t)\right|t^i/i! < +\infty$.
\end{remark}

\section{Discussion}
\label{disc}
We have explored the theory behind the posterior consistency/inconsistency of tail index for Bayesian kernel mixture models, extending  the scope of the vast literature on Bayesian consistency with respect to the weak and strong topology.  We have shown that  examples of inconsistency are extremely common, among the location-scale mixture models with MFM, DPM and NRMI mixture priors. There are special cases in which posterior consistency remains unknown in the DPM and NRMI mixture examples when the marginal base measures of the location and scale parameters meet certain restrictions.


We have also proposed a set of sufficient conditions that lead to posterior tail index consistency, and verified them in a Pareto mixture example. The simple Pareto mixture model is mainly used for illustration, as other heavy tailed kernels with an explicit tail index parameter can also be implemented in a similar manner, such as the inverse gamma kernel, the half Student's $t$ kernel, and the $F$ kernel, although their consistency theory involves extra technical complexity in verifying all those sufficient conditions. It is less obvious to see how models like \eqref{parmix} can be generalized to mixing models with two-sided kernels, since ideally one wants to estimate both the left tail index and the right tail index of a distribution, which can be possibly different. It will be an interesting topic to further study the posterior convergence rates for Model \eqref{parmix} when the true $F_0(x)$ comes from certain nonparametric classes such as the Hall and Welsh class, and compare them with the frequentist adaptive estimators such as \cite{CarKim15} and \cite{BouTho15}, which achieve the minimax rates.


\appendix

\setcounter{equation}{0}
\setcounter{lemma}{0}
\renewcommand{\theequation}{A.\arabic{equation}}
\renewcommand{\thelemma}{A.\arabic{lemma}}

\section{Technical Proofs}


\noindent {\bf Proof of Proposition \ref{borelset}:}\\
\cite{Lan73} Theorem 3.5 has proved that $\mathcal{F}$, the set of all distributions that are absolutely continuous with respect to the Lebesgue measure on $\mathbb{R}$, is a Borel set. Below, we show that the following sets
\begin{align}\label{a2set}
\mathcal{A}_{1}(a) &= \left\{F\in \mathcal{F}:~ \alpha_+(F) \leq a \right\}, \nonumber \\
\mathcal{A}_{2}(a) &= \left\{F\in \mathcal{F}:~ \alpha_+(F) < a \right\},
\end{align}
are Borel sets for any $a\in [0,+\infty]$. First let $a\in [0,+\infty)$. Then for a generic continuous function $g$ on $\mathbb{R}$, we have the following relation
\begin{align}\label{gset}
& \Big\{g:~ \liminf_{x\to+\infty} g(x) \leq a \Big\}  = \Big\{g:~ \sup_{k\in \mathbb{Z}^+,k\geq 2} \inf_{r_j\in \mathbb{Q}^+} g(k+r_j) \leq a \Big\} \nonumber \\
= & \bigcap_{k=2}^{+\infty} \Big\{g:~ \inf_{r_j\in \mathbb{Q}^+} g(k+r_j) \leq a \Big\}  = \left(\bigcup_{k=2}^{+\infty} \left\{g:~ \inf_{r_j\in \mathbb{Q}^+} g(k+r_j) > a \right\}\right)^c \nonumber  \\
= & \left(\bigcup_{k=2}^{+\infty} \bigcup_{q_l\in \mathbb{Q}^+} \left\{g:~ \inf_{r_j\in \mathbb{Q}^+} g(k+r_j) \geq a + q_l \right\}\right)^c \nonumber  \\
= & \left(\bigcup_{k=2}^{+\infty} \bigcup_{q_l\in \mathbb{Q}^+} \bigcap_{r_j\in \mathbb{Q}^+} \left\{g:~ g(k+r_j) \geq a + q_l \right\}\right)^c, \nonumber \\
\text{and }& \Big\{g:~ \liminf_{x\to+\infty} g(x) < a \Big\}  = \Big\{g:~ \sup_{k\in \mathbb{Z}^+,k\geq 2} \inf_{r_j\in \mathbb{Q}^+} g(k+r_j) < a \Big\} \nonumber  \\
= & \bigcap_{k=2}^{+\infty} \Big\{g:~ \inf_{r_j\in \mathbb{Q}^+} g(k+r_j) < a \Big\}  = \left(\bigcup_{k=2}^{+\infty} \left\{g:~ \inf_{r_j\in \mathbb{Q}^+} g(k+r_j) \geq a \right\}\right)^c \nonumber  \\
= & \left(\bigcup_{k=2}^{+\infty} \bigcap_{r_j\in \mathbb{Q}^+} \left\{g:~ g(k+r_j) \geq a \right\}\right)^c,
\end{align}
where $\mathbb{Q}^+$ is the set of all positive rational numbers.

For any $F\in \mathcal{F}$, we have that $-\log \overline F$ is a continuous function on $\mathbb{R}$ (in case $\overline F(x)=0$ for all $x\in [x_1,+\infty)$ with some finite number $x_1$, we can extend the concept of continuity by defining $-\log \overline F(x) = +\infty$ for all $x\geq x_1$). Since $1/\log x$ is continuous for $x\in (1,+\infty)$, we have that the product $-\log \overline F(x) /\log x$ is also continuous on $(1,+\infty)$. For given $x>1,b\geq 0$, we define
\begin{align}\label{dset}
\mathcal{D}(x,b)& = \left\{F \in \mathcal{F}:~ \frac{-\log \overline F(x)}{\log x} \geq b\right\}  = \left\{F\in \mathcal{F}:~ \overline F(x)\leq x^{-b} \right\} \nonumber \\
&= \left\{F\in \mathcal{F}:~ F(x) \geq 1- x^{-b} \right\}
 = \left\{F\in \mathcal{F}:~ F(x) < 1- x^{-b} \right\}^c.
\end{align}
Then \eqref{a2set}, \eqref{gset} and \eqref{dset} together imply that for any $a\in [0,+\infty)$,
\begin{align}\label{a2d2}
\mathcal{A}_{1}(a) & =  \left(\bigcup_{k=2}^{+\infty} \bigcup_{q_l\in \mathbb{Q}^+} \bigcap_{r_j\in \mathbb{Q}^+} \mathcal{D}(k+r_j,a+q_l)\right)^c, \nonumber \\
\mathcal{A}_{2}(a) & =  \left(\bigcup_{k=2}^{+\infty} \bigcap_{r_j\in \mathbb{Q}^+} \mathcal{D}(k+r_j,a)\right)^c.
\end{align}
For any fixed $x>1$ and fixed $p\in (0,1]$, the set $\left\{F:~ F(x)< p\right\}$ is the pre-image of the Borel set $[0,p)$ under the mapping $T_A: F\mapsto F(A)$ for the given Borel set $A=(-\infty, x]$. On the other hand, we know that the Borel sigma-algebra on the space of all distributions $F$ is defined as the smallest sigma-algebra that makes the mapping $F\mapsto F(A)$ measurable for any Borel set $A\subseteq \mathbb{R}$. Using this definition, we know that $\left\{F:~ F(x)< p\right\}$ for fixed $x>1$ and $p\in (0,1]$ is a Borel set. Therefore, $\mathcal{D}(x,b)$ in \eqref{dset} is a Borel set, which further implies that in \eqref{a2d2}, both $\mathcal{A}_{1}(a)$ and $\mathcal{A}_{2}(a)$ are Borel sets for any $a\in [0,+\infty)$.


If $a=+\infty$, then $\mathcal{A}_{1}(+\infty)=\mathcal{F}$ is trivially Borel, and $\mathcal{A}_{2}(+\infty) = \bigcup_{l=2}^{+\infty} \mathcal{A}_{2}(l)$ is also Borel since every $\mathcal{A}_{2}(l)$ is Borel for $l=2,3,\ldots$. Therefore, both $\mathcal{A}_{1}(a)$ and $\mathcal{A}_{2}(a)$ are Borel sets for any $a\in [0,+\infty]$.

Finally we can write $B_{\alpha+}(F,\epsilon) = \mathcal{A}_{2}(\alpha_+(F) +\epsilon ) \cap \left(\mathcal{A}_{1}(\alpha_+(F)-\epsilon)\right)^c$ (in case $\alpha_+(F)-\epsilon<0$, then $\mathcal{A}_{1}(\alpha_+(F)-\epsilon)$ is understood as the empty set). Thus $B_{\alpha+}(F,\epsilon)$ is a Borel set. \hfill $\blacksquare$

\vspace{0.6cm}

In the following, $P_F$ and $E_F$ represent the probability and the expectation under the probability distribution $F$. A random variable $X$ has the decomposition $X=X_+-X_-$, where $X_+=\max(X,0)$ and $X_-=\max(-X,0)$.

\begin{lemma}\label{timom}
(\citealt{ShoWel86}, Theorem 1, Section 7 in Chapter 4) Let $F$ be an univariate distribution on $\mathbb{R}$ with right tail index $\alpha_+(F)$ as defined in \eqref{tidef}. If a random variable $X$ has the cdf $F(x)$, then\\
\begin{align*}
& E_F X_+^m = \begin{cases}
< +\infty & \text{ if } 0< m < \alpha_+(F) \\
= +\infty & \text{ if } m > \alpha_+(F)
\end{cases}
\end{align*}
\end{lemma}

\vspace{.6cm}

\begin{lemma}\label{mnorm}
Let $m>0$.\\
\noindent (i) For any $x,y\in \mathbb{R}$, there exists a constant $C_m$ that only depends on $m$, such that
$$\left[(x + y)_+\right]^m\leq C_m \left(x_+^m + y_+^m\right).$$
\noindent (ii) For any $x\geq 0,y\geq 0$, there exists a constant $c_m$ that only depends on $m$, such that
$$\left(x + y\right)^m\geq c_m \left(x^m + y^m\right).$$
\end{lemma}
\noindent {\bf Proof of Lemma \ref{mnorm}:}\\
(i) For $m\geq 1$, $C_m=2^{m-1}$. For $m\in (0,1)$, $C_m=1$.\\
(ii) Let $f(t) = t^m+(1-t)^m$ and $t\in [0,1]$. If $m\geq1$, then $\max_{t\in [0,1]}f(t)=1$ and set $c_m=1$. If $m\in (0,1)$, then $\max_{t\in [0,1]}f(t)=2^{1-m}$ and set $c_m=2^{m-1}$. Now let $t=x/(x+y)$ and the conclusion follows. \hfill$\blacksquare$
\vspace{.6cm}

\begin{lemma}\label{mixbound}
Suppose $f$ is a density drawn from Model \eqref{model2} with cdf $F$. Let $K(\cdot)$ be the cdf of $k(\cdot)$. Then
\begin{align}
& E_F X_+^m \geq  c_m \left(E_{G_{\mu}}\mu_+^m \cdot \overline K(0) + E_K X_+^m \cdot E_{G_{\mu,\sigma}}\left[\sigma^{m}I(\mu\geq 0)\right]\right),  \label{mbound1} \\
& E_F X_+^m \geq  C_m^{-1} E_K X_+^m \cdot E_{G_{\sigma}}\sigma^{m} - E_{G_{\mu}}\mu_-^m, \label{mbound2} \\
& E_F X_+^m \leq  C_m \left(E_{G_{\mu}}\mu_+^m + E_K X_+^m \cdot E_{G_{\sigma}}\sigma^{m}\right),  \label{mbound3}
\end{align}
where $\overline K(0) = P_K(X\geq 0)$ (the probability of $X\geq 0$ if $X$ has the density $k(x)$), $G_{\mu}$ and $G_{\sigma}$ are the marginal distributions of $G_{\mu,\sigma}$, and $c_m, C_m$ are defined in Lemma \ref{mnorm}.
\end{lemma}

\noindent {\bf Proof of Lemma \ref{mixbound}:}\\
Let $I(\cdot)$ be the indicator function. We have
\begin{align}\label{mk1}
E_F X_+^m & = \iint x^mI(x\geq 0) \frac{1}{\sigma}k\left(\frac{x-\mu}{\sigma}\right) \ud G(\mu,\sigma) \ud x \nonumber \\
& = \iint (\mu+\sigma y)^m I(\mu+\sigma y\geq 0) k(y) \ud G(\mu,\sigma) \ud y.
\end{align}
Then we give lower and upper bounds for \eqref{mk1}. Notice that $\sigma\geq 0$ always holds and $I(\mu+\sigma y\geq 0)  \geq I(\mu\geq 0) I(y\geq 0)$. Based on \eqref{mk1} and part (ii) of Lemma \ref{mnorm}, we have:
\begin{align*}
&E_F X_+^m \geq  \iint (\mu+\sigma y)^m I(\mu\geq 0)I(y\geq 0) k(y) \ud G(\mu,\sigma) \ud y \\
&\geq \iint c_m(\mu_+^m + \sigma^m y_+^m)I(\mu\geq 0)I(y\geq 0) k(y) \ud G(\mu,\sigma) \ud y\\
&=  c_m \left(E_{G_{\mu}}\mu_+^m \cdot \overline K(0)+ E_K X_+^m \cdot E_{G_{\mu,\sigma}}\left[\sigma^{m}I(\mu\geq 0)\right]\right),
\end{align*}
which is \eqref{mbound1}.

On the other hand, since $(-\mu)_+=\mu_-$, part (i) of Lemma \ref{mnorm} implies
\begin{align*}
& (\sigma y)_+^m \leq C_m\left[(\mu+\sigma y)_+^m + (-\mu)_+^m\right] \\
\implies & (\mu+\sigma y)_+^m \geq C_m^{-1} \sigma^my_+^m - \mu_-^m.
\end{align*}
This together with \eqref{mk1} gives
\begin{align*}
&E_F X_+^m = \iint (\mu+\sigma y)_+^m  k(y) \ud G(\mu,\sigma) \ud y \\
&\geq \iint \left[C_m^{-1} \sigma^my_+^m - \mu_-^m\right] k(y) \ud G(\mu,\sigma) \ud y \\
&\geq C_m^{-1} E_K X_+^m \cdot E_{G_{\sigma}}\sigma^{m} - E_{G_{\mu}}\mu_-^m,
\end{align*}
which is \eqref{mbound2}.

By part (i) of Lemma \ref{mnorm}
\begin{align*}
&E_F X_+^m =  \iint [(\mu+\sigma y)_+]^m  \ud G(\mu,\sigma) \ud y \\
&\leq \iint C_m(\mu_+^m + \sigma^m y_+^m) k(y) \ud G(\mu,\sigma) \ud y  =  C_m \left(E_{G_{\mu}}\mu_+^m + E_K X_+^m \cdot E_{G_{\sigma}}\sigma^{m}\right),
\end{align*}
which is \eqref{mbound3}. \hfill$\blacksquare$
\vspace{.6cm}


\noindent {\bf Proof of Theorem \ref{mfmcounter}:}\\
For Model \eqref{model3}, the marginal distributions $G_{\mu}$ and $G_{\sigma}$ are both finite mixtures at the points $\mu_{i=1}^N$ and $\sigma_{i=1}^N$ respectively. Because $G_{0}(\mu,\sigma)$ is a continuous distribution, we have $0\leq E_{G_{\mu}}\mu_+^m < +\infty$, $0\leq E_{G_{\mu}}\mu_-^m < +\infty$ and $0<E_{G_{\sigma}}\sigma^m <+\infty$ for all $m>0$. We can use Lemma \ref{timom} to determine the relation between $\alpha_+(F)$ and $\alpha_+(K)$. According to Lemma \ref{mixbound}, whether $E_F X_+^m$ is finite or not for a given $m$ is solely determined by whether $E_K X_+^m$ is finite or not. The analysis goes as follows: \\
\noindent (i) If $\alpha_+(K)=+\infty$, then by Lemma \ref{timom} $E_K X_+^m < +\infty$ for all $m>0$. The upper bound \eqref{mbound3} implies that $E_F X_+^m<+\infty$ for all $m>0$. Hence, $\alpha_+(F)=+\infty$ by Lemma \ref{timom}. \\
\noindent (ii) If $\alpha_+(K)=0$, then by Lemma \ref{timom}, $E_K X_+^m =+\infty$ for all $m>0$. The lower bound \eqref{mbound2} implies that $E_F X_+^m=+\infty$ for all $m>0$. Then by setting $m=0$ in (ii) of Lemma \ref{timom} we can see that $\alpha_+(F)=0$. \\
\noindent (iii) If $\alpha_+(K) \in (0,+\infty)$, then $E_K X_+^m<+\infty$ for $m<\alpha_+(K)$ and $E_K X_+^m=+\infty$ for $m>\alpha_+(K)$. Then by \eqref{mbound3}, $E_F X_+^m <+\infty$ for $m<\alpha_+(K)$, and by \eqref{mbound2}, $E_K X_+^m=+\infty$ for $m>\alpha_+(K)$. Apply Lemma \ref{timom} and we can see that $\alpha_+(F)=\alpha_+(K)$. \\
In sum, $\alpha_+(F)=\alpha_+(K)$ in all three cases and thus $\alpha_+(F)$ is almost surely a singleton. \hfill $\blacksquare$

\vspace{.6cm}

The homogenous NRMI with L\'evy intensity $\rho(\ud v)$ and base measure $H_0$ defined in Section 3 can be expressed as $H(x) =\sum_{i\geq 1} J_i\delta_{X_i}(x)$ where $J_i=\widetilde J_i/ \sum_{i\geq 1}\widetilde J_i$. Equivalently, the cdf $H(x)$ also has the representation $H(x)=S(H_0(x))/S(1)$, where $\{S(t),t\geq 0\}$ is a subordinator with L\'evy intensity measure $\rho(\ud v)$ (see for example \citealt{Regetal03}). As a result, the function $\Psi$ defined in Theorem \eqref{nrmcounter1} is the Laplace exponent of the subordinator $S(t)$. The conditions $\int_0^{\infty}\rho(\ud v)=+\infty$ and $\int_0^{\infty}(1-e^{-v})\rho(\ud v)<+\infty$ guarantees that $0<S(1)<+\infty$ almost surely.

The following proposition is a combination of Theorem 1 in \cite{Fri67} and Lemmas 4 and 5 in \cite{FriPru71}.

%

\begin{prop}\label{Friprop}
(\citealt{Fri67, FriPru71}) Suppose $\{S(t),t\geq 0\}$ is a subordinator with L\'evy intensity measure $\rho(\ud v)$ for $v\in \mathbb{R}^+$.  Define the functionals $R_L(h)=\liminf_{t\to 0+} S(t)/h(t)$ and $R_U(h)=\limsup_{t\to 0+} S(t)/h(t)$. \\
	\noindent (i) For $\gamma>0$, let $h_{\gamma}(x)$ be the same as defined in \eqref{hgamma}.
	Then
	\begin{align*}
	& R_L(h_{\gamma}) \leq \gamma \text{ a.s. if }\gamma<1,\\
	& R_L(h_{\gamma}) \geq \gamma-1 \text{ a.s. if }\gamma>1.
	\end{align*}
	\noindent (ii) If $h(x)$ is locally convex in $x\in [0,\epsilon)$ for some small $\epsilon>0$, then
	\begin{align*}
	& R_U(h) = 0 \text{ a.s. if }\int_0^\epsilon \rho[h(x),+\infty)\ud x < +\infty,\\
	& R_U(h) = +\infty \text{ a.s. if }\int_0^\epsilon \rho[h(x),+\infty)\ud x = +\infty.
	\end{align*}
\end{prop}

\vspace{.6cm}

\noindent {\bf Proof of Theorem \ref{nrmcounter1}:} \\
The proof is a direct application of Proposition \ref{Friprop}. \\
\noindent (i) By the stationary increment property of subordinators, $S(1-t)$ has the same distribution as $S(1)-S(t)$ for $t\in (0,1)$. Therefore for $\gamma>1$ and $h_{\gamma}$ defined in \eqref{hgamma}, part (i) of Proposition \ref{Friprop} implies
$$\liminf_{t\to 1-} \frac{S(1)-S(t)}{h_{\gamma}(1-t)}\geq \gamma-1 \text{ a.s.}$$
Let $t= H_0(x)$ and we have
\begin{align}\label{gam1}
&\liminf_{x\to+\infty} \frac{\overline H(x)S(1)}{h_{\gamma}\left(\overline H_0(x)\right)}\geq \gamma-1 \text{ a.s.}
\end{align}
since $H(x) = S(H_0(x))/S(1)$. Our assumptions $\int_0^{\infty}\rho(\ud v)=+\infty$ and $\int_0^{\infty}(1-e^{-v})\rho(\ud v)<+\infty$ guarantee that $0<S(1)<+\infty$ almost surely. Therefore, we conclude from \eqref{gam1} that almost surely for all such NRMI $H$,
\begin{align*}
&\liminf_{x\to+\infty} \frac{\overline H(x)}{h_{\gamma}\left(\overline H_0(x)\right)}>0 \text{ a.s.}
\end{align*}
As $x\to+\infty$, the function $\overline H(x)/h_{\gamma}\left(\overline H_0(x)\right)$ is almost surely lower bounded by a positive constant, which implies that the function $\log\left[\overline H(x)/h_{\gamma}\left(\overline H_0(x)\right)\right]$ is almost surely lower bounded by a finite constant. Hence it follows that
\begin{align}\label{gam2}
&\liminf_{x\to+\infty} \frac{\log\left[\overline H(x)/h_{\gamma}\left(\overline H_0(x)\right)\right]}{\log x} \geq 0 \text{ a.s.}
\end{align}
For the right tail index of $H$, we can use \eqref{gam2} and the condition on $h_{\gamma}(\cdot)$ to obtain that
\begin{align*}
\alpha_+(H)&=\liminf_{x\to+\infty}\frac{-\log \overline H(x)}{\log x} \\
& = \liminf_{x\to+\infty}\left\{\frac{-\log h_{\gamma}\left(\overline H_0(x)\right)}{\log x} +
 \frac{\log \left[h_{\gamma}\left(\overline H_0(x)\right)/\overline H(x) \right] }{\log x}  \right\} \\
& \leq \liminf_{x\to+\infty}\frac{-\log h_{\gamma}(\overline H_0(x))}{\log x} - \liminf_{x\to +\infty}\frac{\log\left[\overline H(x)/h_{\gamma}\left(\overline H_0(x)\right)\right]}{\log x} \\
& \leq \liminf_{x\to+\infty}\frac{-\log h_{\gamma}(\overline H_0(x))}{\log x} = 0.
\end{align*}
Therefore $\alpha_+(H)=0$.\\
\noindent (ii) For such $h(x)$ that satisfies (a)(b)(c), by similar argument as above, we apply part (ii) of Proposition \ref{Friprop} and obtain that
$$\limsup_{x\to+\infty} \frac{\overline H(x)S(1)}{h\left(\overline H_0(x)\right)}=0 \text{ a.s.}$$
which implies that almost surely for all such NRMI $H$,
\begin{align*}
\limsup_{x\to+\infty} \frac{\overline H(x)}{h\left(\overline H_0(x)\right)}=0 \text{ a.s.}
\end{align*}
Therefore, we have
$$\liminf_{x\to+\infty} \log \frac{h\left(\overline H_0(x)\right)}{\overline H(x)} = +\infty \text{ a.s.}$$
and hence
$$\liminf_{x\to+\infty} \frac{\log \left[h\left(\overline H_0(x)\right)/\overline H(x)\right]}{\log x} \geq 0 \text{ a.s.}$$
We finally combine this with the condition (c) and conclude that
\begin{align*}
\alpha_+(H)&=\liminf_{x\to+\infty}\frac{-\log \overline H(x)}{\log x} \\
&= \liminf_{x\to+\infty} \left\{\frac{-\log h\left(\overline H_0(x)\right)}{\log x} + \frac{\log\left[h\left(\overline H_0(x)\right)/\overline H(x)\right]}{\log x}\right\} \\
&\geq \liminf_{x\to+\infty} \frac{-\log h\left(\overline H_0(x)\right)}{\log x} + \liminf_{x\to+\infty} \frac{\log\left[h\left(\overline H_0(x)\right)/\overline H(x)\right]}{\log x} \\
&\geq \liminf_{x\to+\infty} \frac{-\log h\left(\overline H_0(x)\right)}{\log x} = +\infty.
\end{align*}
which means $\alpha_+(H)=+\infty$. \hfill$\blacksquare$

\vspace{.6cm}

\noindent {\bf Proof of Theorem \ref{nrmcounter2}:} \\
First we note that because $G_0(\mu,\sigma)$ is a continuous probability measure, if $G(\cdot,\cdot)$ is a homogenous NRMI with L\'evy intensity $\rho(\ud v)G_0(\ud \mu,\ud \sigma)$, then using the stick-breaking representation, we have that the two marginal distributions $G_{\mu}$ and $G_{\sigma}$ are also homogenous NRMIs with L\'evy intensities $\rho(\ud v)G_{0,\mu}(\ud \mu)$ and $\rho(\ud v)G_{0,\sigma}(\sigma)$ respectively. Given the conclusion of Theorem \ref{nrmcounter1}, we have that if $G_{0,\mu}$ or $G_{0,\sigma}$ satisfies (i) of Theorem \ref{nrmcounter1}, then $\alpha_+(G_{\mu})=0$ or $\alpha_+(G_{\sigma})=0$; if $G_{0,\mu}$ or $G_{0,\sigma}$ satisfies (ii) of Theorem \ref{nrmcounter1}, then $\alpha_+(G_{\mu})=+\infty$ or $\alpha_+(G_{\sigma})=+\infty$.

Since $k(\cdot)$ has part of the support in $\mathbb{R}^+$, $E_K X_+^m>0$ for any $m>0$. We will examine the existence of moments $E_F X_+^m$ with $F$ from Model \ref{model2} for any $m>0$, and use Lemma \ref{timom} to determine $\alpha_+(F)$. Similar to the proof of Theorem \ref{mfmcounter}, we can analysis $E_F X_+^m$ using the lower bounds and the upper bound from Lemma \ref{mixbound}.\\

\noindent (i) If $\alpha_+(G_{\mu})=0$, then $E_{G_{\mu}}\mu_+^m=+\infty$ for all $m>0$. Also note that $\overline K(0)>0$ and $E_K X_+^m >0$ since $k(\cdot)$ has full support in $\mathbb{R}$. Therefore, by the lower bound \eqref{mbound1}, $E_F X_+^m=+\infty$ for all $m>0$ since $\overline K(0)>0$, $E_K X_+^m >0$, and $E_{G_{\mu,\sigma}}\left[\sigma^{m}I(\mu\geq 0)\right]\geq 0$. This implies $\alpha_+(F)=0$ by Lemma \ref{timom}. \\

\noindent (ii) If $\alpha_+(G_{\mu})=+\infty$ and $\alpha_+(G_{\sigma})=0$, then for all $m>0$, $E_{G_{\mu}} \mu_+^m <+\infty$ and $E_{G_{\sigma}} \sigma^m = +\infty$. Because we have assumed that $G_{0,\mu}$ is symmetric about zero, this implies that $E_{G_{\mu}} \mu_-^m <+\infty$ for all $m>0$. Also $E_K X_+^m >0$ for all $m>0$. Therefore by the lower bound \eqref{mbound2}, $E_F X_+^m=+\infty$. This again implies $\alpha_+(F)=0$ by Lemma \ref{timom}. \\

\noindent (iii) If $\alpha_+(G_{\mu})=+\infty$ and $\alpha_+(G_{\sigma})=+\infty$, then for all $m>0$, $E_{G_{\mu}} \mu_+^m <+\infty$ and $E_{G_{\sigma}} \sigma^m< +\infty$. This can be further separated into three scenarios: (a) $\alpha_+(K)\in (0,+\infty)$, then if $m\in (0,+\infty)$ and $m <\alpha_+(K)$, $E_KX_+^m <+\infty$ and $E_F X_+^m<+\infty$ by the upper bound \eqref{mbound3}; if $m\in (0,+\infty)$ and $m >\alpha_+(K)$, $E_KX_+^m = +\infty$ and $E_F X_+^m=+\infty$ by the lower bound \eqref{mbound2}. Hence $\alpha_+(F) =\alpha_+(K)$ by Lemma \ref{timom}. (b) $\alpha_+(K)=0$, then $E_K X_+^m = +\infty$ for all $m\in (0,+\infty)$ and $E_F X_+^m = +\infty$ for all $m\in (0,+\infty)$ by the lower bound in \eqref{mbound2}. (c) $\alpha_+(K)=+\infty$, then $E_K X_+^m < +\infty$ for all $m\in (0,+\infty)$ and $E_F X_+^m < +\infty$ for all $m\in (0,+\infty)$ by the upper bound in \eqref{mbound3}. We conclude that in all three scenarios $\alpha_+(F)=\alpha_+(K)$.\\

The results from different scenarios can be summarized as
$$\alpha_+(F)=\min\left\{\alpha_+(G_{\mu}),\alpha_+(G_{\sigma}),\alpha_+(K)\right\},$$
which is always a fixed number. Therefore, $\alpha_+(F)$ is almost surely a singleton, if both $G_{0,\mu}(\mu)$ and $G_{0,\sigma}(\sigma)$ satisfy either (i) or (ii) in Theorem \ref{nrmcounter1}. \hfill$\blacksquare$

\vspace{.6cm}

\noindent {\bf Proof of Theorem \ref{dpcounter}:}\\
We will show that for a measure $H_0$ on $\mathbb{R}$ \\
\noindent (a) If $\limsup_{x\to +\infty}\overline H_0(x)\cdot \left(\log x /\log\log x\right) = +\infty$, then part (i) of Theorem \ref{nrmcounter1} holds;\\
\noindent (b) If $\limsup_{x\to +\infty}\overline H_0(x) \cdot \left[(\log x)\cdot (\log\log x)^{\delta}\right] = 0$ for some $\delta>1$, then part (ii) of Theorem \ref{nrmcounter1} holds.\\
If both $\overline G_{0,\mu}$ and $\overline G_{0,\sigma}$ satisfy either (a) or (b), i.e. we can replace $H_0$ with $\overline G_{0,\mu}$ and $\overline G_{0,\sigma}$, then the right tail indices of $G_{\mu}$ and $G_{\sigma}$ are either $0$ or $+\infty$, and the conclusion of Theorem \ref{dpcounter} follows directly from Theorem \ref{nrmcounter2}.\\
To show (a) and (b), we use the similar arguments as in \cite{DosSel82}. We note that a cdf $H(x)$ on $\mathbb{R}$ drawn from $\DP(a,H_0)$ can be written as a normalized Gamma process with L\'evy intensity $\rho(\ud v)H_0(\ud x) = a v^{-1} e^{-v} \ud v H_0(\ud x)$. The Laplace exponent for $\rho$ is $\Psi(s)=a\log(1+s)$ and its inversion is $\Psi^{-1}(u)= e^{u/a}-1$. Thus for any given $\gamma>1$, the function \eqref{hgamma} in Proposition \ref{Friprop} is given by
$$h_{\gamma}(x)= \frac{\log|\log x|}{\exp\left(\frac{\gamma \log |\log x|}{ax}\right)-1}.$$
We have $\lim_{x\to 0+}h_{\gamma}(x) = 0$, and $h_{\gamma}(x)\in (0,1/2)$ for $x\in[0,\epsilon)$ for small enough $\epsilon>0$.

Now by the condition $\limsup_{x\to +\infty}\overline H_0(x)/(\log\log x/\log x) = +\infty$, there exists a positive sequence $x_j$ that increases to $+\infty$ as $j\to +\infty$, such that for any $C>2$, $1/16 > \overline H_0(x_j) > C \log\log x_j/\log x_j$ and $x_j>\exp(C^2)$ as long as $j>J(C)$ for some large integer $J(C)$. Therefore for all $j>J(C)$,

\begin{align*}
& \frac{-\log h_{\gamma}\left(\overline H_0(x_j)\right)}{\log x_j} =  \frac{\log\left[\exp\left(\frac{\gamma \log |\log \overline H_0(x_j)|}{a\overline H_0(x_j)}\right)-1\right] - \log \log \left|\log \overline H_0(x_j)\right|}{\log x_j} \\
&\leq \frac{\gamma\log \left|\log \overline H_0(x_j)\right|}{a\overline H_0(x_j)\log x_j}
\leq \frac{\gamma \log \left|\log \log x_j - \log C - \log\log\log x_j\right|}{a C\log\log x_j} \leq \frac{\gamma\log \log\log x_j}{aC\log \log x_j}.
\end{align*}
As $j\to +\infty$, this upper bound converges to 0. Together with the fact that $h_{\gamma}(x)\in (0,1/2)$ for $x\in[0,\epsilon)$, we obtain that $\liminf_{x\to+\infty}-\log h_{\gamma}\left(\overline H_0(x)\right)/\log x =0$. This is exactly the condition in part (i) of Theorem \ref{nrmcounter1}. Thus (a) is proved. \\

For (b), we set $h(x) = \exp\left[-\left(x|\log x|^{\delta'}\right)^{-1}\right]$ for some $1<\delta'<\delta$. This function is convex in $[0,\epsilon)$ for small enough $\epsilon>0$. It also satisfies $\lim_{x\to 0+}h(x) = 0$, and $h(x)\in (0,1/2)$ for $x\in[0,\epsilon)$ for small enough $\epsilon>0$. Due to the lower and upper bounds $ae^{-1}\log{1/u}\leq \rho[u,+\infty)\leq a\log (1/u)+ae^{-1}$ (see \citealt{DosSel82}) and $\delta'>1$, we have $\int_0^{\epsilon} \rho[h(x),+\infty)\ud x<+\infty$. Furthermore, if $\limsup_{x\to +\infty}\overline H_0(x) \cdot \left[(\log x)\cdot (\log\log x)^{\delta}\right]=0$, then for any $C>2$ and all sufficiently large $x$, $\overline H_0(x) < \min\left(1/\left[C(\log x)\cdot (\log\log x)^{\delta}\right],1/2\right)$. Therefore for sufficiently large $x$,
\begin{align*}
& \frac{-\log h(\overline H_0(x))}{\log x} =  \frac{1}{\overline H_0(x) \left|\log \overline H_0(x)\right|^{\delta'} \log x}
\geq \frac{C (\log x)\cdot (\log\log x)^{\delta}}{\left\{\log \left[C\log x \cdot (\log\log x)^\delta\right]\right\}^{\delta'}\log x} \\
& \geq \frac{C (\log x) (\log\log x)^{\delta}}{\left(2\log \log x\right)^{\delta'}\log x} = \frac{C}{2^{\delta'}}(\log\log x)^{\delta-\delta'}\to +\infty,
\end{align*}
which implies that $\liminf_{x\to+\infty} -\log h(\overline H_0(x))/\log x = +\infty$. This is exactly the condition in part (ii) of Theorem \ref{nrmcounter1}. Thus (b) is proved. \hfill$\blacksquare$
\vspace{.6cm}

\noindent {\bf Proof of Theorem \ref{nggpcounter}:}\\
We will show that for a measure $H_0$ on $\mathbb{R}$ \\
\noindent (a) If $\limsup_{x\to +\infty}\overline H_0(x) \cdot x^{\delta} = +\infty$ for all $\delta>0$, then part (i) of Theorem \ref{nrmcounter1} holds;\\
\noindent (b) If $\limsup_{x\to +\infty}\overline H_0(x) \cdot  x^{\delta} =0 $ for all $\delta>0$, then part (ii) of Theorem \ref{nrmcounter1} holds.\\
If both $\overline G_{0,\mu}$ and $\overline G_{0,\sigma}$ satisfy either (a) or (b), i.e. we can replace $H_0$ with $\overline G_{0,\mu}$ and $\overline G_{0,\sigma}$, then the right tail indices of $G_{\mu}$ and $G_{\sigma}$ are either $0$ or $+\infty$, and the conclusion of Theorem \ref{nggpcounter} follows directly from Theorem \ref{nrmcounter2}.\\
To show (a), we note that for the L\'evy process with intensity $\rho(\ud v) = \frac{a}{\Gamma(1-\kappa)}v^{-\kappa-1}e^{-\tau v}\ud v$, its Laplace exponent is $\Psi(s)=\frac{a}{\kappa}\left[(s+\tau)^\kappa-\tau^\kappa\right]$, and its inverse is $\Psi^{-1}(u)=$\\
 $\left[\kappa u/a+\tau^\kappa\right]^{1/\kappa}-\tau$. Thus for any given $\gamma>1$, the function \eqref{hgamma} is given by
$$h_{\gamma}(x)= \frac{\log |\log x|}{\left[\frac{\kappa \log |\log x|}{ax} +\tau^\kappa\right]^{1/\kappa}-\tau}$$
We have $\lim_{x\to 0+}h_{\gamma}(x) = 0$, and $h_{\gamma}(x)\in (0,1/2)$ for $x\in[0,\epsilon)$ for small enough $\epsilon>0$.

Now by the condition $\limsup_{x\to +\infty}\overline H_0(x)\cdot x^{\delta} = +\infty$ for all $\delta>0$, we have the following conclusion: for any given $\delta>0$, there exists a positive sequence $x_j$ that increases to $+\infty$ as $j\to +\infty$, such that $x_j>16$ and $\min\left(\kappa\log\log 16/(a\tau^\kappa), 1/16\right) > \overline H_0(x_j) > x_j^{-\delta}$ as long as $j>J$ for some large integer $J$. Such choice of $x_j$ guarantees that $\log\log(\delta\log x_j)>\log\log\left|\log \overline H_0(x_j) \right|> \log\log\log 16 > 0$, and
\begin{align*}
&\frac{\kappa\log \left|\log\overline H_0(x_j)\right|}{a\overline H_0(x_j)} >
\frac{\kappa\log \log 16}{a\overline H_0(x_j)} > \tau^\kappa.
\end{align*}
Therefore for all $j>J$,
\begin{align*}
& \frac{-\log h_{\gamma}(\overline H_0(x_j))}{\log x_j} \leq  \frac{\kappa^{-1} \log\left[ \frac{2\kappa \log \left|\log\overline H_0(x_j)\right|}{a\overline H_0(x_j)}\right]- \log \log \left|\log \overline H_0(x_j)\right|}{\log x_j} \\
&= \frac{-\kappa^{-1} \log \overline H_0(x_j) + \left(\kappa^{-1}-1\right)\log \log \left|\log \overline H_0(x_j)\right| + \kappa^{-1}\log (2\kappa/a)}{\log x_j}\\
&\leq \kappa^{-1}\delta + \left(\kappa^{-1}-1\right)\frac{\log\log\left(\delta\log x_j\right)}{\log x_j} + \frac{\log (2\kappa/a)}{\kappa\log x_j}.
\end{align*}
In the last display, the second and the third terms converge to zero as $j\to+\infty$. The first term can be made arbitrarily small if $\delta$ is made small. Therefore, we have shown that $\liminf_{x\to+\infty} -\log h_{\gamma}(\overline H_0(x))/\log x=0$. Thus (a) is proved.\\

For (b), we have the following bound for $\rho[u,+\infty)$:
\begin{align*}
\rho[u,+\infty) & = \frac{a}{\Gamma(1-\kappa)}\int_u^1 t^{-\kappa-1}e^{-\tau t}\ud t +  \frac{a}{\Gamma(1-\kappa)}\int_1^\infty t^{-\kappa-1}e^{-\tau t}\ud t \\
& \leq \frac{a}{\Gamma(1-\kappa)}\int_u^1 t^{-\kappa-1}\ud t + \frac{a}{\Gamma(1-\kappa)}\int_1^\infty e^{-\tau t}\ud t \\
& =  \frac{a}{\kappa\Gamma(1-\kappa)} \left(\frac{1}{u^{\kappa}}-1\right) + \frac{ae^{-\tau}}{\tau\Gamma(1-\kappa)}.
\end{align*}
Therefore if $\int_0^\epsilon \frac{1}{h(x)^{\kappa}} \ud x <+\infty$, then $\int_0^\epsilon \rho[h(x),+\infty) \ud x <+\infty$. Let $h(x) = (x|\log x|^{\eta})^{1/\kappa}$ for some $\eta>1$. For $\kappa\in (0,1)$, this $h(x)$ is convex and increasing in $[0,\epsilon)$ and satisfies $\lim_{x\to 0+}h(x) = 0$, and $h(x)\in (0,1/2)$ for $x\in[0,\epsilon)$ if $\epsilon$ is sufficiently small.

On the other hand, if $\limsup_{x\to +\infty}\overline H_0(x) \cdot  x^{\delta} =0 $ for all $\delta>0$, then $\overline H_0(x) \leq x^{-\delta}$ for any given $\delta$ for all sufficiently large $x$, which implies that
\begin{align*}
& \frac{-\log h(\overline H_0(x))}{\log x} \geq \frac{-\log h(x^{-\delta})}{\log x}
= \frac{\delta}{\kappa} - \frac{\eta \log \left(\delta \log x\right) }{\kappa \log x}.
\end{align*}
On the right-hand side of the last display, the second term converges to zero as $x\to+\infty$. The first term can be made arbitrarily large if $\delta$ is made large. Therefore we have shown that $\liminf_{x\to+\infty} -\log h(\overline H_0(x))/\log x =+\infty$. Thus (b) is proved. \hfill$\blacksquare$




\vspace{.6cm}

\noindent {\bf Proof of Theorem \ref{utest1}:}\\
Let $s_n$ be a positive sequence such that $B_n \prec s_n \prec \overline\alpha_n^{-1} \log n$, whose existence is guaranteed by Condition (iv). For $\epsilon>0$, we define the test $\Phi_n = I\left(|\hat\alpha_{s_n}-\alpha_{0+}|\geq \epsilon/2\right)$ with $\hat\alpha_{s_n}$ given by \eqref{alphatest}. Let $p_{s_n}=P_F(X>e^{s_n})$ be the population mean of $\hat p_{s_n}$ and let $\alpha_{s_n} = \log(p_{s_n}) - \log(p_{s_n+1})$. Note that $p_{s_n}$ and $\alpha_{s_n}$ implicitly depend on $F$. We complete the proof in two steps.

\noindent {\bf Step 1}: Show $E_{F_0}\Phi_n \to 0$ as $n\to \infty$.

We have
\begin{align}\label{ute1}
E_{F_0} \Phi_n & = P_{F_0} \left(|\hat \alpha_{s_n}-\alpha_{0+}|\geq \frac{\epsilon}{2} \right) \nonumber \\
& \leq P_{F_0} \left(|\hat \alpha_{s_n}-\alpha_{s_n}|\geq \frac{\epsilon}{4} \right) +   P_{F_0} \left(| \alpha_{s_n}-\alpha_{0+}|\geq \frac{\epsilon}{4} \right).
\end{align}
The first term in \eqref{ute1} can be bounded by Lemma 2 and equation (4.2) of \cite{CarKim15}:
\begin{align}\label{ute11}
& P_{F_0}\Big(|\hat \alpha_{s_n} - \alpha_{s_n}| \geq \frac{\epsilon}{4}\Big) \leq 2\exp\left(-\frac{np_{s_n+1}\epsilon^2}{576}\right),
\end{align}
where $p_{s_n+1}=P_{F_0}(X>e^{s_n+1}) = e^{-\alpha_{0+}(s_n+1)}L_0(e^{s_n+1})$. Since $L_0$ is slowly varying, as $n\to \infty$, eventually $L_0(e^{s_n+1})\geq e^{-\delta(s_n+1)}$ for arbitrarily small $\delta>0$. By Condition (iv) of Theorem \ref{utest1},  $s_n \prec \log n/\alpha_{0+}$ and hence $np_{s_n+1}\geq \exp\left(\log n-(\alpha_{0+}+\delta)(s_n+1)\right) \to +\infty$, which implies that the righthand side of \eqref{ute11} goes to zero as $n\to \infty$.

The second term in \eqref{ute1} is not stochastic. We have
\begin{align*}
|\alpha_{s_n}-\alpha_{0+}|& = \left|\log(p_{s_n}) - \log(p_{s_n+1}) -\alpha_{0+}\right| \\
& = \left|\log \frac{L_0(e^{s_n})}{L_0(e^{s_n+1})}\right| \to 0,
\end{align*}
because $L_0$ is slowly varying and $s_n\to \infty$. Therefore both terms on the righthand side of \eqref{ute1} converge to zero as $n\to \infty$.

\noindent {\bf Step 2}: Show $\sup_{F\in B^c_{\alpha+}(F_0,\epsilon)\cap \mathcal{F}_n} E_{F}(1-\Phi_n)\to 0$ as $n\to \infty$, where we let $\mathcal{F}_n = \mathcal{F}_{1n}\cap \mathcal{F}_{2n}\cap\mathcal{F}_{3n}$. By Conditions (i)-(iii), it is clear that $\Pi_n(\mathcal{F}_n^c)\leq \Pi_n(\mathcal{F}_{1n}^c)+ \Pi_n(\mathcal{F}_{2n}^c) + \Pi_n(\mathcal{F}_{3n}^c) \leq e^{-c_1n}+e^{-c_2n}+e^{-c_3n} \leq e^{-c'n}$ where $c'=\min(c_1,c_2,c_3)/2$.

For every $F\in B^c_{\alpha+}(F_0,\epsilon)\cap\mathcal{F}_n$, we have $|\alpha_+(F)-\alpha_{0+}|>\epsilon$. Therefore
\begin{align}\label{ute2}
&  E_F(1-\Phi_n) = P_F \left(|\hat \alpha_{s_n} - \alpha_{0+}|\leq \frac{\epsilon}{2} \right) \nonumber \\
& = P_F \left(|\hat \alpha_{s_n} - \alpha_{0+}|\leq \frac{\epsilon}{2}, |\hat \alpha_{s_n} - \alpha_+(F)|< \frac{\epsilon}{2} \right) \nonumber \\
& ~~ +P_F \left(|\hat \alpha_{s_n} - \alpha_{0+}|\leq \frac{\epsilon}{2}, |\hat \alpha_{s_n} - \alpha_+(F)|\geq \frac{\epsilon}{2} \right)\nonumber \\
& \leq P_F \left(| \alpha_+(F)- \alpha_{0+}|< \epsilon\right) + P_F \left( |\hat \alpha_{s_n} - \alpha_+(F)|\geq \frac{\epsilon}{2} \right) \nonumber \\
& =P_F \left( |\hat \alpha_{s_n} - \alpha_+(F)|\geq \frac{\epsilon}{2} \right) \leq P_F \left( |\hat \alpha_{s_n} - \alpha_{s_n}|\geq \frac{\epsilon}{4} \right) + P_F \left( |\alpha_{s_n} - \alpha_+(F)|\geq \frac{\epsilon}{4} \right).
\end{align}
We only need to show that both terms on the righthand side of \eqref{ute2} converge to zero uniformly over all $F\in B^c_{\alpha+}(F_0,\epsilon)\cap\mathcal{F}_n$ as $n\to \infty$. For a fixed $F$, the first term can be bounded by Lemma 2 and equation (4.2) of \cite{CarKim15} again as
$$P_{F}\left(|\hat \alpha_{s_n} - \alpha_{s_n}| \geq \frac{\epsilon}{4}\right) \leq 2\exp\left(-\frac{np_{s_n+1}\epsilon^2}{576}\right).$$
To obtain uniform convergence for the righthand side, we only need the quantity $np_{s_n+1}$ to be uniformly bounded below for all $F\in B^c_{\alpha+}(F_0,\epsilon)\cap\mathcal{F}_n$. Using Conditions (i)-(iii), we can obtain the following uniform lower bound:
\begin{align*}
& np_{s_n+1} = n e^{-\alpha_+(F) (s_n+1)} L_F(e^{s_n+1}) =  n e^{-\alpha_+(F) (s_n+1)} L_F(x_0) \exp\left(\int_{x_0}^{e^{s_n+1}} \frac{h_F(t)}{t}\ud t\right) \\
& \geq \exp\left(\log n - \overline \alpha_n (s_n+1) - c_L\log n - \int_{x_0}^{e^{s_n+1}} \frac{B_n}{t(\log t)^{1+\tau_n}} \ud t \right) \\
& = \exp\left((1-c_L)\log n - \overline \alpha_n (s_n+1) - \frac{B_n}{\tau_n (\log x_0)^{\tau_n}}+\frac{B_n}{\tau_n (s_n+1)^{\tau_n}}  \right) \\
&\geq \exp\left((1-c_L)\log n - \overline \alpha_n (s_n+1)  - \frac{B_n}{\tau_n}\right),
\end{align*}
where we use $x_0\geq e$ and hence $\log x_0 \geq 1 $ in the last inequality. Condition (i) says $1-c_L>0$. By our choice of $s_n$, we have $\log n \succ  \overline \alpha_n (s_n+1)$, and Condition (iv) implies $\log n \succ  B_n/\tau_n$. Therefore we have obtained that uniformly over all $F\in B^c_{\alpha+}(F_0,\epsilon)\cap\mathcal{F}_n$, $P_{F}\left(|\hat \alpha_{s_n} - \alpha_{s_n}| \geq \frac{\epsilon}{4}\right)$ converges to zero as $n\to\infty$.

For the second term in \eqref{ute2}, we have
\begin{align*}
& |\alpha_{s_n}-\alpha_+(F)| = \left|\log p_{s_n} - \log p_{s_n+1} - \alpha_+(F) \right| \\
={}& \left|\log e^{-\alpha_+(F) s_n}L_F(e^{s_n}) - \log e^{-\alpha_+(F) (s_n+1)}L_F(e^{s_n+1})  - \alpha_+(F) \right| \\
={}& \left|\log L_F(e^{s_n}) - \log L_F(e^{s_n+1})  \right|  = \left| \int_{e^{s_n}}^{e^{s_n+1}} \frac{h_F(x)}{x} \ud x  \right| \\
\leq{}&  \int_{e^{s_n}}^{e^{s_n+1}} \frac{\overline h_n(x)}{x} \ud x =  \int_{e^{s_n}}^{e^{s_n+1}} \frac{B_n}{x(\log x)^{1+\tau_n}} \ud x \\
={}& \frac{B_n}{\tau_n s_n^{\tau_n}} \left[1-\left(\frac{s_n}{1+s_n}\right)^{\tau_n}\right]
=\frac{B_n}{\tau_n s_n^{\tau_n}}  \left[1-\exp\left(-\tau_n\log\frac{1+s_n}{s_n}\right)\right] \\
\leq{}&  \frac{B_n}{\tau_n s_n^{\tau_n}} \cdot \tau_n\log\left(1+\frac{1}{s_n}\right) \leq \frac{B_n}{s_n^{1+\tau_n}},
\end{align*}
where we have used $1-e^{-t}\leq t$ for $t>0$ and $\log(1+t)\leq t$ for $t>0$. Since $1\preceq B_n\prec s_n$, we have $B_n/s_n^{1+\tau_n} \to 0$ as $n\to\infty$. Therefore the probability $P_F \left( |\alpha_{s_n} - \alpha_+(F)|\geq \frac{\epsilon}{4} \right)$ is zero for all large $n$ uniformly over all $F\in B^c_{\alpha+}(F_0,\epsilon)\cap\mathcal{F}_n$.  \hfill$\blacksquare$

\vspace{.6cm}

\begin{lemma}\label{f0exp}
A distribution $F \in  \mathcal{CM}_e\cap \mathcal{P}_2$ if and only if the survival function $\overline F$ and the density $f$ take the form
\begin{align*}
& \overline F(x) = w x^{-\alpha} + (1-w) \int_{\alpha (1+\beta)}^{+\infty} x^{-u} \ud H (u), \\
& f(x) = w \alpha x^{-(\alpha+1)} + (1-w) \int_{\alpha (1+\beta)}^{+\infty} u x^{-(u+1)} \ud H (u),
\end{align*}
where $w\in (0,1]$, $\alpha>0$, $\beta>0$, $H$ is a probability measure whose support is in $[\alpha(1+\beta),+\infty)$.
\end{lemma}

\noindent {\bf Proof of Lemma \ref{f0exp}:} \\
We only need to show the expression of $\overline F(x)$ since the expression of $f(x)$ will follow directly. We note that if a distribution $F$ has a density $f$, then $F(e^t)$ for $t\in [0,+\infty)$ is also a cdf and $\overline F(e^t)$ is a survival function. By the Hausdorff-Bernstein-Widder theorem, $\overline F(e^t)$ is a completely monotone function on $t\in[0,+\infty)$ if and only if it is the Laplace transformation of some probability distribution $G$ on $(0,+\infty)$, i.e.
$$\overline F(e^t) =\int_0^{\infty} e^{-u t}\ud G(u),$$
which is equivalent to say that for $x\in [1,+\infty)$,
\begin{align}\label{h31}
&\overline F(x) = \int_0^{\infty} x^{-u}\ud G(u).
\end{align}
Therefore to prove the conclusion of the lemma, it only remains to show that $F \in \mathcal{P}_2$ if and only if the probability measure $G$ has the decomposition
\begin{align}\label{h32}
&G(u) = w\delta_{\alpha} + (1-w) H(u),
\end{align}
for some $w\in (0,1]$, $\alpha>0$ and some probability measure $H$ supported on $[\alpha(1+\beta),+\infty)$.

It is clear that if $G$ has the form \eqref{h32} (so that $\overline F$ has the form in the lemma), then $F \in \mathcal{P}_2$. Conversely, if $F$ satisfies \eqref{h31} and meanwhile $\overline F(x)=Cx^{-\alpha}+O(x^{-\alpha(1+\beta)})$ for some $\alpha,\beta>0$ and $C>0$, then we show that:\\
\noindent (i) $G(0,\alpha)=0$;\\
\noindent (ii) $G(\alpha,\alpha(1+\beta))=0$.

If (i) does not hold, then there exists a set $I\subset (0,\alpha)$ such that $G(I)>0$. Thus by Fatou's lemma,
$$\liminf_{x\to+\infty} x^{\alpha}\int_{I} x^{-u}\ud G(u)
\geq \int_{I} \liminf_{x\to+\infty} x^{\alpha-u} \ud G(u) = +\infty,$$
which contradicts the fact that
$$\limsup_{x\to+\infty}x^\alpha \int_{I} x^{-u}\ud G(u) \leq \limsup_{x\to+\infty}x^\alpha \int_0^\infty x^{-u}\ud G(u)
\leq \limsup_{x\to+\infty} x^\alpha \overline F(x) = C  <+\infty.$$
Hence (i) must hold true.

Suppose the point mass of $G$ at $\alpha$ is $w$ ($w\geq 0$). Then since $G$ is a probability measure, $G-w\delta_{\alpha}$ is a nonnegative measure supported on $(\alpha,+\infty)$. For any set $I\subseteq (\alpha,+\infty)$, by Fatou's lemma,
$$\limsup_{x\to +\infty} x^\alpha \int_I x^{-u} \ud \left[G(u)-w\delta_{\alpha}(u)\right]
\leq  \int_I \limsup_{x\to +\infty}x^{\alpha-u} \ud \left[G(u)-w\delta_{\alpha}(u)\right]=0.$$
Now we compare this with the format of $\overline F(x)=Cx^{-\alpha}+O(x^{-\alpha(1+\beta)})$. Because $G$ is a probability measure, we must have $w=C\in (0,1]$ and for $I=(\alpha,+\infty)$, \\
$\int_I  x^{-u} \ud \left[G(u)-w\delta_{\alpha}(u)\right] = O(x^{-\alpha(1+\beta)})$, which is equivalent to
\begin{align}\label{limsup2}
&\limsup_{x\to+\infty}x^{\alpha(1+\beta)}\int_I x^{-u}\ud \left[G(u)-C\delta_{\alpha}(u)\right] < +\infty.
\end{align}
If (ii) does not hold, then there exists another set $I'\subseteq (\alpha,\alpha(1+\beta))$ such that $G(I')>0$. Then by Fatou's lemma,
$$\liminf_{x\to+\infty} x^{\alpha(1+\beta)}\int_{I'} x^{-u}\ud \left[G(u)-w\delta_{\alpha}(u)\right]
\geq \int_{I'} \liminf_{x\to+\infty} x^{\alpha(1+\beta)-u} \ud \left[G(u)-w\delta_{\alpha}(u)\right] = +\infty,$$
which contradicts \eqref{limsup2}. Thus we have shown that if $F \in \mathcal{CM}_e\cap \mathcal{P}_2$, then both (i) and (ii) have to hold and $G$ satisfies \eqref{h32}. \hfill $\blacksquare$
\vspace{.6cm}


\noindent {\bf Proof of Theorem \ref{cmkl}:}\\
We prove the theorem in a similar way to the proof of the exponential mixture model in Theorem 16 of \cite{WuGho08}.

By Lemma \ref{f0exp}, we assume that the true density function has the form
\begin{align}\label{f0form}
& f_0(x) = w_0 \alpha_0 x^{-(\alpha_0+1)} + (1-w_0) \int_{\alpha_0(1+\beta_0)}^{\infty} \alpha x^{-(\alpha+1)} \ud G_0 (\alpha),
\end{align}
where $\alpha_0$ is short for $\alpha_{0+}(F)$ and $G_0$ is supported on $[\alpha_0(1+\beta_0),+\infty)$. Without causing confusion, we also denote a generic $f$ from $\mathcal{CM}_e\cap \mathcal{P}_2$ by $f_{w_1,\alpha_1,H}$ which takes the form of Model (4.3) in the main paper.

We use $\overline \Pi_n$ to denote the overall prior measure, including the prior on $w_1,\alpha_1,H_1$.


The KL condition for Model (4.3) in the main paper is satisfied if for any $\epsilon>0$, there exists sets $\mathcal{W}\subset (0,1]$ (for $w_1$), $\mathcal{A}\subset (0,+\infty)$ (for $\alpha_1$) and $\mathcal{H}$ (for $H$) that do not depend on $n$, such that
\begin{align}\label{liminf3set}
&G_{w}(\mathcal{W})>0,~ G_{\alpha}(\mathcal{A})>0, ~ \Pi(\mathcal{H};\bxi,H_0)>0,
\end{align}
and meanwhile for all $(w_1,\alpha_1,H)\in \mathcal{W}\times \mathcal{A}\times \mathcal{H}$ and all sufficiently large $n$,
\begin{align}\label{kltotal}
\int_1^{+\infty} f_0(x) \log \frac{f_0(x)}{f_{w_1,\alpha_1,H}(x)} \ud x < \epsilon.
\end{align}
This is because as $n\to\infty$, $\underline w_n\to 0$, $\overline\alpha_n\to+\infty$, and
\begin{align*}
& \liminf_{n\to \infty} \overline \Pi_n(\mathcal{W}\times \mathcal{A}\times \mathcal{H})= \liminf_{n\to \infty} \left(G_{w}\cdot I_{[\underline w_n,1]}\right)(\mathcal{W})\times \left(G_{\alpha}\cdot I_{[0,\overline\alpha_n]}\right)(\mathcal{A}) \times \Pi(\mathcal{H};\bxi,H_0) \\
&=\liminf_{n\to \infty} \frac{G_{w}(\mathcal{W})}{1-G_w([0,\underline w_n))}\times \liminf_{n\to \infty} \frac{G_{\alpha}(\mathcal{A})}{1-G_{\alpha}((\overline\alpha_n,+\infty))} \times \Pi(\mathcal{H};\bxi,H_0) \\
&= G_{w}(\mathcal{W})\times G_{\alpha}(\mathcal{A})\times \Pi(\mathcal{H};\bxi,H_0)>0.
\end{align*}

To show \eqref{liminf3set} and \eqref{kltotal}, by Theorem 1 of \cite{WuGho08}, we prove the following 3 relations to obtain the conclusion for all $\epsilon\in (0,1/2)$ and all sufficiently large $n$:\\
\noindent (i) $\int_1^{\infty} f_0(x) \log \frac{f_0(x)}{f_{w_0,\alpha_0,G_1}(x)} \ud x < \epsilon/3$ for a distribution $G_1$ supported on $[\alpha_0(1+\beta_0),a]$ with a fixed $a>\alpha_0(1+\beta_0)$;  \\
\noindent (ii) $\int_1^{\infty} f_0(x) \log \frac{f_{w_0,\alpha_0,G_1}(x)}{f_{w_1,\alpha_1,G_1}(x)} \ud x < \epsilon/3$ for all $(w_1,\alpha_1)\in \mathcal{W}\times \mathcal{A}$ such that
$G_{\alpha}(\mathcal{A})>0$ and $ G_w(\mathcal{W})>0$; \\
\noindent (iii) $\int_1^{\infty} f_0(x) \log \frac{f_{w_1,\alpha_1,G_1}(x)}{f_{w_1,\alpha_1,H}(x)} \ud x < \epsilon/3$ for all $H'\in \mathcal{H}'$ and all $(w_1,\alpha_1)\in \mathcal{W}\times \mathcal{A}$, and \\
$\Pi(\mathcal{H};\bxi,H_0) >0$.
\vspace{.3cm}

\noindent Check (i): Let $a>\alpha_0(1+\beta_0)+3$ whose value will be chosen later. For each fixed $a$, there exists a large integer $n(a)$ such that for all large $n\geq n(a)$, $a < \overline \alpha_n$. Let $G_1(A)=G_0(A)/G_0([\alpha_0(1+\beta_0), a])$ for any set $A\subseteq [\alpha_0(1+\beta_0),a]$. Then for every $x\geq 1$, $f_{w_0,\alpha_0,G_1}(x)$ converges pointwise to $f_0(x)$ by taking $a\to +\infty$ and then $n \geq n(a)\to \infty$. For sufficiently large $a$, we pick a fixed number $\alpha_2\in (\alpha_0(1+\beta_0),a)$ such that $G_0([\alpha_0(1+\beta_0),\alpha_2])\geq 1/2$. We notice that for fixed $x>1$, the function $\alpha x^{-(\alpha+1)}$ attains its maximum when $\alpha=1/\log x$, increases on $(0,1/\log x]$ and decreases on $(1/\log x,+\infty)$. Based on this property, the following relations hold:
\begin{align*}
& \frac{1}{2} \alpha_2 x^{-(\alpha_2+1)} \leq \int_{\alpha_0(1+\beta_0)}^{a} ux^{-(u+1)}\ud G_1(u)\leq \alpha_0(1+\beta_0) x^{-(\alpha_0(1+\beta_0)+1)}  \\
& \qquad \text{ if } x\geq e^{1/[\alpha_0(1+\beta_0)]}, \\
& \frac{1}{2} \min\left(\alpha_0(1+\beta_0) x^{-(\alpha_0(1+\beta_0)+1)},\alpha_2 x^{-(\alpha_2+1)}\right) \leq \int_{\alpha_0(1+\beta_0)}^{a} ux^{-(u+1)}\ud G_1(u) \leq \frac{1}{ex\log x} \\
& \qquad \text{ if } e^{1/\alpha_2}\leq x< e^{1/[\alpha_0(1+\beta_0)]}, \\
& \frac{1}{2} \alpha_0(1+\beta_0) x^{-(\alpha_0(1+\beta_0)+1)}  \leq \int_{\alpha_0(1+\beta_0)}^{a} ux^{-(u+1)}\ud G_1(u) \leq \alpha_2 x^{-(\alpha_2+1)} \\
& \qquad \text{ if } 1\leq x< e^{1/\alpha_2}.
\end{align*}
This gives a lower bound and an upper bound for $f_{w_0,\alpha_0,G_1}(x)$:
\begin{align*}
& f_{w_0,\alpha_0,G_1}(x) \geq w_0\alpha_0 x^{-(\alpha_0+1)} + \frac{1-w_0}{2} \min\left(\alpha_0(1+\beta_0) x^{-(\alpha_0(1+\beta_0)+1)},\alpha_2 x^{-(\alpha_2+1)}\right) \\
&:=g_1(x), \\
& f_{w_0,\alpha_0,G_1}(x) \leq w_0\alpha_0 x^{-(\alpha_0+1)} + (1-w_0)\Big[\alpha_0(1+\beta_0) x^{-(\alpha_0(1+\beta_0)+1)} + \frac{1}{ex\log x} \\
& + \alpha_2 x^{-(\alpha_2+1)}\Big] \leq \frac{1}{ex\log x} + w_0\alpha_0 + (1-w_0)[\alpha_0(1+\beta_0)+\alpha_2] := g_2(x),
\end{align*}
where $g_1(x)$ and $g_2(x)$ are defined to be the lower and upper bounds of $f_{w_0,\alpha_0,G_1}(x)$ as above.
Therefore, we can obtain an upper bound for $|\log f_{w_0,\alpha_0,G_1}(x)|$ for all $x\in (1,+\infty)$:
\begin{align}\label{g12bound}
|\log f_{w_0,\alpha_0,G_1}(x)| \leq & |\log g_1(x)| + |\log g_2(x)|.
\end{align}
Since clearly both $\log x$ and $|\log \log x|$ are $f_0$-integrable, $|\log g_1(x)|$ and $|\log g_2(x)|$ are also $f_0$-integrable, and so is $|\log f_{w_0,\alpha_0,G_1}(x)|$. By the dominated convergence theorem, as $a\to +\infty$ and $n \geq n(a)\to \infty$, we have $\int_1^{\infty} f_0(x) \log \frac{f_0(x)}{f_{w_0,\alpha_0,G_1}(x)} \ud x \to 0$. Therefore, for any given $\epsilon>0$, we can pick a fixed $a$ now, such that for all $n\geq n(a)$, $\int_1^{\infty} f_0(x) \log \frac{f_0(x)}{f_{w_0,\alpha_0,G_1}(x)} \ud x <\epsilon/3$, which is the conclusion of Part (i). Note that $G_1$ is now a fixed distribution and does not depend on $n$.
\vspace{.3cm}

\noindent Check (ii): We show (ii) for the $G_1$ constructed in Part (i). Let $\mathcal{W}=[w_0-\eta,w_0+\eta]$ and $\mathcal{A} = [\alpha_0-\eta,\alpha_0+\eta]$ for some $\eta\in (0,1)$. Then since $f_{w_1,\alpha_1,G_1}(x)$ is a continuous function of $(w_1,\alpha_1)$ at $(w_0,\alpha_0)$, $f_{w_1,\alpha_1,G_1}(x)\to f_{w_0,\alpha_0,G_1}(x)$ pointwise in $x$, uniformly for all $w_1\in \mathcal{W}$ and $\alpha_1\in \mathcal{A}$ as $\eta\to 0$. Hence there exists $\eta>0$ such that $f_{w_0,\alpha_0,G_1}(x)/2\leq f_{w_1,\alpha_1,G_1}(x)\leq 2f_{w_0,\alpha_0,G_1}(x)$ for all $w_1\in \mathcal{W}$ and all $\alpha_1\in \mathcal{A}$. Since $|\log f_{w_0,\alpha_0,G_1}(x)|$ is $f_0$-integrable, it implies that $|\log f_{w_1,\alpha_1,G_1}(x)| \leq |\log [f_{w_0,\alpha_0,G_1}(x)/2]| + |\log [2f_{w_0,\alpha_0,G_1}(x)]|$ is also $f_0$-integrable. Therefore by \eqref{g12bound} and the dominated convergence theorem, as $\eta\to 0$, $\int_1^{\infty} f_0(x) \log \frac{f_{w_0,\alpha_0,G_1}(x)}{f_{w_1,\alpha_1,G_1}(x)} \ud x \to 0$. Therefore, for any given $\epsilon>0$, we can choose a fixed constant $\eta\in (0,1)$, such that $\int_1^{\infty} f_0(x) \log \frac{f_{w_0,\alpha_0,G_1}(x)}{f_{w_1,\alpha_1,G_1}(x)} \ud x < \epsilon/3$ for all $(w_1,\alpha_1)\in \mathcal{W}\times \mathcal{A}$, which is the conclusion of Part (ii). Since in Model (4.3) in the main paper, the support of $G_{\alpha}$ will include $\mathcal{A}=[\alpha_0-\eta,\alpha_0+\eta]$ and the support of $G_{w}$ will include $\mathcal{W}=[w_0-\eta,w_0+\eta]$ as $n\to\infty$, we have that $G_{\alpha}(\mathcal{A})>0$ and $G_w(\mathcal{W})>0$.
\vspace{.3cm}

\noindent Check (iii): The argument is similar to the proof of Lemma 3 in \cite{WuGho08}. We split the integral in Part (iii) into two parts:
\begin{align}\label{Itotal}
& \int_1^{\infty} f_0(x) \log \frac{f_{w_1,\alpha_1,G_1}(x)}{f_{w_1,\alpha_1,H}(x)} \ud x \nonumber \\
&\leq \int_1^{C_1} f_0(x) \log \frac{f_{w_1,\alpha_1,G_1}(x)}{f_{w_1,\alpha_1,H}(x)} \ud x +\int_{C_1}^{\infty} f_0(x) \log \frac{f_{w_1,\alpha_1,G_1}(x)}{f_{w_1,\alpha_1,H}(x)} \ud x \nonumber \\
& := I_1+I_2.
\end{align}
We bound $I_1$ and $I_2$ separately. For all $H\in \mathcal{H}'_1$ and all $(w_1,\alpha_1)\in \mathcal{W}\times\mathcal{A}$, $I_2$ can be upper bounded as
\begin{align*}
I_2& \leq \int_{C_1}^{\infty} f_0(x) \log \frac{f_{w_1,\alpha_1,G_1}(x)}{w_1\alpha_1x^{-(\alpha_1+1)}} \ud x \\
&\leq \int_{C_1}^{\infty} f_0(x) \left|\log f_{w_1,\alpha_1,G_1}(x) \right| + \int_{C_1}^{\infty} f_0(x) \left[|\log(w_1\alpha_1)|+(\alpha_1+1)\log x\right]\ud x.
\end{align*}
Clearly both $|\log f_{w_1,\alpha_1,G_1}(x)|$ and $|\log(w_1\alpha_1)|+(\alpha_1+1)\log x$ are $f_0$-integrable
uniformly for all $(w_1,\alpha_1)\in \mathcal{W}\times\mathcal{A}$. Therefore we can choose $C_1$ sufficiently large, such that $I_2< \epsilon/6$.


To bound $I_1$, we notice that
\begin{align}\label{I11}
I_1 &\leq \sup_{x\in [1,C_1]} \left|\frac{f_{w_1,\alpha_1,G_1}(x)}{f_{w_1,\alpha_1,H}(x)}-1\right|\leq \sup_{x\in [1,C_1]}\left|\frac{\int_{\alpha_1}^{+\infty}k(x;u) \ud [G_1(u)-H(u)]}{f_{w_1,\alpha_1,H}(x)}\right| \nonumber \\
& \leq \frac{\sup_{x\in [1,C_1]}\left|\int_{\alpha_1}^{+\infty} k(x;u) \ud [G_1(u)-H(u)]\right| }
{\inf_{x\in [1,C_1]} f_{w_1,\alpha_1,H}(x)}.
\end{align}
For any distribution $H$, all $(w_1,\alpha_1)\in \mathcal{W}\times\mathcal{A}$,
\begin{align}\label{I12}
& \inf_{x\in [1,C_1]} f_{w_1,\alpha_1,H}(x) \geq \inf_{x\in [1,C_1]} w_1\alpha_1 C_1^{-(\alpha_1+1)} \nonumber\\
& \geq (w_0-\eta)(\alpha_0-\eta)C_1^{-(\alpha_0+\eta+1)}:=C_2.
\end{align}
Let $\mathcal{H}_1=\{H_1\sim \Pi(H;\bxi,H_0): H_1([0,a-\alpha_0-1])>1- C_2\epsilon/(24a)\}$. Then since $\Pi$ is a homogeneous NRMI, $\Pi(\mathcal{H}_1;\bxi,H_0)>0$. Let $\mathcal{H}'_1= \{H:H(\alpha) = H_1(\alpha-\alpha_1),H_1\in \mathcal{H}_1\}$. Let $\mathcal{D}=[\alpha_0-\eta,a]$. Then $H(\mathcal{D})\geq H_1([0,a-\alpha_0-\eta])\geq H_1([0, a-\alpha_0-1])>1- C_2\epsilon/(24a)$ for any $H\in \mathcal{H}'_1$. Because $G_1(\mathcal{D})=1>1-C_2\epsilon/(24a)$, we know that $\mathcal{H}'_1$ is an open neighborhood of $G_1$. Then for all $H\in \mathcal{H}'_1$,
\begin{align}\label{I13}
\sup_{x\in [1,C_1]} \left|\int_{\mathcal{D}^c} k(x;u) \ud [G_1(u)-H(u)]\right|  &\leq a\left(G_1(\mathcal{D}^c)+H(\mathcal{D}^c)\right) < \frac{C_2\epsilon}{24}
\end{align}
since $G_1(\mathcal{D}^c)=0$. Because the kernel $k(x;\alpha)$ is equicontinuous on $\mathcal{D}$, by the Arzela-Ascoli theorem, there exist $N$ points $x_1,\ldots,x_N\in [1,C_1]$ such that for any $x\in [1,C_1]$, $\sup_{u\in \mathcal{D}}|k(x;u)-k(x_i;u)|< C_2\epsilon/24$ for some $i=1,\ldots,N$. Now we choose a smaller open neighborhood $\mathcal{H}'_2\subseteq \mathcal{H}'_1$ for $H$ such that $\max_{i=1,\ldots,N}\left|\int_{\mathcal{D}} k(x_i;u)\ud [G_1(u)-H(u)]\right|<C_2\epsilon/24$ for all $H\in \mathcal{H}'_2$. The correspondingly open set for $H_1$ is $\mathcal{H}_2=\{H_1: H_1(\alpha)=H(\alpha+\alpha_1), H\in \mathcal{H}'_2\}$), and it satisfies $\Pi(\mathcal{H}_2; \bxi,H_0)>0$ since $\mathcal{D}$ is a fixed interval and $N$ is finite. Then for any $x\in [1,C_1]$ and $H\in \mathcal{H}'_2$, there exists some $x_i$ ($i=1,\ldots,N$) such that
\begin{align}\label{I14}
&\left|\int_{\alpha_1}^{+\infty} k(x;u) \ud [G_1(u)-H(u)]\right| \nonumber \\
\leq{}& \left|\int_{\mathcal{D}} k(x;u) \ud [G_1(u)-H(u)]\right| + \left|\int_{\mathcal{D}^c} k(x;u) \ud [G_1(u)-H(u)]\right| \nonumber \\
\leq{}& \left|\int_{\mathcal{D}} [k(x;u)-k(x_i;u)] \ud G_1(u)\right| + \left|\int_{\mathcal{D}} [k(x;u)-k(x_i;u)] \ud H(u)\right| \nonumber \\
& + \left|\int_{\mathcal{D}} k(x_i;u) \ud [G_1(u)-H(u)]\right| + \left|\int_{\mathcal{D}^c} k(x;u) \ud [G_1(u)-H(u)]\right| \nonumber \\
<{}& \frac{C_2\epsilon}{8} + \left|\int_{\mathcal{D}^c} k(x;u) \ud [G_1(u)-H(u)]\right|.
\end{align}
For $\epsilon\in (0,1/2)$, \eqref{I13} and \eqref{I14} together give us
\begin{align}\label{I15}
& \sup_{x\in [1,C_1]} \left|\int_{\alpha_1}^{+\infty} k(x;u) \ud [G_1(u)-H(u)]\right| \leq \frac{C_2\epsilon}{6},
\end{align}
for all $H\in \mathcal{H}'$. We combine \eqref{I11}, \eqref{I12} and \eqref{I15} to obtain that $I_1<(C_2\epsilon)/(6C_2)=\epsilon/6$. Therefore we have shown that in \eqref{Itotal}, for all $(w_1,\alpha_1)\in \mathcal{W}\times\mathcal{A}$ and all $H\in \mathcal{H}_3$,
$$\int_1^{\infty} f_0(x) \log \frac{f_{w_1,\alpha_1,G_1}(x)}{f_{w_1,\alpha_1,H}(x)} \ud x\leq I_1+I_2< \frac{\epsilon}{3}.$$
Finally we set $\mathcal{H}=\mathcal{H}_2$ such that $\Pi(\mathcal{H}; \bxi,H_0)>0$. Hence (iii) is proved. (i)-(iii) together imply \eqref{liminf3set} and \eqref{kltotal}, which further implies the conclusion of Theorem \ref{cmkl}. \hfill $\blacksquare$

\vspace{.6cm}

\noindent {\bf Proof of Theorem \ref{parthm}:}\\
Condition (i) implies (KL) by Theorem \ref{cmkl}. We first check the condition (PT) for Model (4.3) in the main paper. For any $x\geq e$ and any distribution $F$ drawn from Model (4.3) in the main paper, we have
\begin{align*}
\limsup_{x\to+\infty} \frac{-\log \overline F(x)}{\log x} & = \limsup_{x\to+\infty} \frac{-\log \left[w_1x^{-\alpha_1} + (1-w_1)\int_{\alpha_1}^{\infty} x^{-\alpha} dH(\alpha)\right]}{\log x} \\
&\leq \limsup_{x\to+\infty} \frac{-\log \left(w_1x^{-\alpha_1} \right)}{\log x} = \limsup_{x\to+\infty} \frac{\alpha_1\log x -\log w_1}{\log x} = \alpha_1, \\
\liminf_{x\to+\infty} \frac{-\log \overline F(x)}{\log x} & = \liminf_{x\to+\infty} \frac{-\log \left[w_1x^{-\alpha_1} + (1-w_1)\int_{\alpha_1}^{\infty} x^{-\alpha} dH(\alpha)\right]}{\log x} \\
&\geq \liminf_{x\to+\infty} \frac{-\log \left[w_1x^{-\alpha_1}+ (1-w_1)x^{-\alpha_1}\right]}{\log x} = \liminf_{x\to+\infty} \frac{\alpha_1\log x}{\log x} = \alpha_1.
\end{align*}
Therefore, the limit exists and $\lim_{x\to+\infty} -\log \overline F(x)/\log x=\alpha_1$, which means that (PT) holds with $\alpha_+(F)=\lim_{x\to+\infty} -\log \overline F(x)/\log x=\alpha_1$.

Now we only need to show (UT) by checking the conditions of Theorem \ref{utest1}. For a distribution $F$ drawn from Model (4.3) in the main paper, the slowly varying function $L_F$ can be written as
\begin{align*}
&L_F(x) = w_1 + (1-w_1)\int_{\alpha_1}^{+\infty} x^{\alpha_1-\alpha} \ud H(\alpha).
\end{align*}
To check Condition (i) of Theorem \ref{utest1}, we notice that $L_F(x_0)\geq w_1 \geq \underline w_n$, which satisfies $\underline w_n \geq n^{-c_L}$ for all sufficiently large $n$ given Condition (iii) of Theorem \ref{parthm}. This relation holds for arbitrary $x_0\geq e$. Therefore we will choose the exact value $x_0$ later when we check Condition (ii) of Theorem \ref{utest1}.

Because $\alpha_+(F)=\alpha_1$ and $\alpha_1$ is drawn from $G_{\alpha}$ with support $(0,\overline \alpha_n]$, Condition (iii) of Theorem \ref{utest1} is satisfied by the same assumption on $\overline \alpha_n$ in Condition (iii) of Theorem \ref{parthm}.

Next we check the uniform bound on the function $h_F$, in Condition (ii) of Theorem \ref{utest1}, for the cases of DP and NGGP in Condition (ii) (a) and (b) of Theorem \ref{parthm}, respectively. We first notice that
\begin{align}\label{hfbound1}
\left|h_F(x)\right|& = \left|\frac{xL'_F(x)}{L_F(x)}\right|=\left|\frac{(1-w_1)\int_{\alpha_1}^{+\infty} (\alpha_1-\alpha)x^{\alpha_1-\alpha} \ud H(\alpha)}{w_1+ (1-w_1) \int_{\alpha_1}^{+\infty} x^{\alpha_1-\alpha} \ud H(\alpha)}\right| \leq w_1^{-1} \int_0^{\infty} ux^{-u} \ud H_1(u),
\end{align}
where $H_1$ is as specified in Model (4.3) in the main paper such that $H_1(u) = H(u+\alpha_1)$ for $u>0$.

\vspace{.3cm}

\noindent (a) Suppose that $\Pi$ is $\DP(a,H_0)$ and the base measure $H_0$ satisfies the conditions in Condition (ii)(a). Recall that the L\'evy intensity of this Dirichlet process is $\rho(\ud v)=av^{-1}e^{-v} \ud v$ and it satisfies $\rho[u,+\infty) \leq a\log (1/u) + ae^{-1}$ for $u\in (0,1)$ (see \citealt{DosSel82}). The distribution $H_1$ from $\DP(a,H_0)$ can be defined as $H_1(u)=S(H_0(u))/S(1)$ for all $u\in \mathbb{R}$ where $S(t)$ for $t\in [0,1]$ is the subordinator with the L\'evy intensity $\rho(\ud v)$.

Now let $h(x)=\exp\left[-x^{-1/(1+\delta)}\right]$ for some $\delta\in (0,d_1)$, where $d_1$ is from Condition (ii)(a). Then $\lim_{x\to 0+}h(x)=0$, $h(x)$ is convex and $h(x)\in (0,1/2)$ for $x\in [0,\epsilon_1]$ for a small enough $\epsilon_1\in (0,1/2)$. Furthermore, $\int_0^{\epsilon} \rho[h(x),+\infty) \ud x <+\infty$. Therefore, according to Part (ii) of Proposition A.1, $\limsup_{t\to 0+} S(t)/h(t)=0$ almost surely. Since $H_0$ has no point mass at zero, $\lim_{u\to 0+} H_0(u)=0$ and hence $\limsup_{u\to 0+} S(H_0(u))/h(H_0(u))=0$ almost surely. Equivalently, since $0<S(1)<+\infty$ almost surely, this implies that $\limsup_{u\to 0+} H_1(u)/h(H_0(u))=0$ almost surely. Hence by Condition (ii)(a), there exists a small constant $0<\epsilon_2\leq \min(c_1,\exp[-(D_1/\epsilon_1)^{1/(1+d_1)}])$ such that $H_1(u)\leq h(H_0(u))$ almost surely for all $u \in (0,\epsilon_2]$, and meanwhile $H_0(u)\leq D_1[\log(1/u)]^{-(1+d_1)}$\\ $\leq D_1 [\log(1/\epsilon_2)]^{-(1+d_1)}\leq \epsilon_1$ for all $u \in (0,\epsilon_2]$.

Now for some small $\eta\in (0,\epsilon_2]$ whose value will be chosen later (which will depend on $x$) and for all $x\geq e$, we bound the last integral in \eqref{hfbound1} as
\begin{align}\label{boundbypart}
& \int_0^{\infty} u x^{-u} \ud H_1(u)  =  \int_0^{\eta} ux^{-u} \ud H_1(u) + \int_{\eta}^{\infty} ux^{-u} \ud H_1(u) \nonumber \\
&\stackrel{(i)}{\leq} \eta \int_0^{\eta} x^{-u} \ud H_1(u) + \int_{\eta}^{\infty} (x/2)^{-u} \ud H_1(u) \nonumber \\
&\stackrel{(ii)}{\leq} \eta x^{-u}  H_1(u) \Big|_0^{\eta} + \eta(\log x) \int_0^{\eta} x^{-u} H_1 (u) \ud u + (x/2)^{-\eta} \nonumber \\
&= \eta x^{-\eta}H_1(\eta) + \eta (\log x) \int_0^{\eta} x^{-u} H_{1}(u) \ud u  + (x/2)^{-\eta} \nonumber \\
&\stackrel{(iii)}{\leq} \eta(\log x) \int_0^{\eta} x^{-u} H_{1}(u) \ud u  + 2(x/2)^{-\eta},
\end{align}
where in (i) we used the fact that $0<u\leq \eta$ in the first term and $u< 2^u$ for all $u>0$ in the second term; in (ii) we applied the integration by parts for the integral on $(0,\eta)$, used the fact that $H_1$ has no point mass at $u=0$ because the base measure $H_0$ has no point mass at $u=0$, and used $x/2\geq e/2>1$; in (iii) we used the fact that $\eta <2^{\eta}$ and $H_1(\eta)\leq 1$. We bound the two terms in \eqref{boundbypart} separately. For the first term, we have that for all $\eta\in (0,\epsilon_2]$ and any $x\geq e$,
\begin{align}\label{term10}
& (\log x) \int_0^{\eta} x^{-u} H_1(u) \ud u \leq (\log x) \int_0^{\eta} x^{-u} h(H_0(u)) \ud u \leq (\log x) \int_0^{\eta} x^{-u} h(H_0(\eta)) \ud u \nonumber \\
& \leq h\left(H_0(\eta)\right) \int_0^{\eta} (\log x)e^{-u(\log x)}  \ud u \leq  h\left(H_0(\eta)\right) \int_0^{\infty} (\log x)e^{-u(\log x)}  \ud u =  h\left(H_0(\eta)\right), \text{ a.s.}
\end{align}
where we used the fact that $h$ is increasing on $(0,\epsilon_1]$ and $H_0(u)$ is non-decreasing in $u$. Under Condition (i)(a), for all $\eta\in (0,\epsilon_2]$ and any $x\geq e$, \eqref{term10} implies that
\begin{align}\label{dphterm11}
& (\log x) \int_0^{\eta} x^{-u} H_1(u) \ud u \leq h\left(D_1 \left[\log (1/\eta)\right]^{-(1+d_1)}\right) \nonumber \\
&=\exp\left\{-D_1^{-1/(1+\delta)}[\log(1/\eta)]^{(1+d_1)/(1+\delta)}\right\}  \text{ a.s.}
\end{align}
Now we choose $\eta = \exp\left\{-\left[D_1^{1/(1+d_1)}(2\log\log x)^{(1+\delta)/(1+d_1)}\right]\right\}$ such that $\eta\leq \epsilon_2$ for all $x\geq x_{1}$. This holds with the constant $x_1=\exp\{\exp[D_1^{-1/(1+\delta)}(\log(1/\epsilon_2))^{(1+d_1)/(1+\delta)}/2]\}$. With this choice of $\eta$, for all $x\geq \max(e,x_1)$, \eqref{dphterm11} implies that
\begin{align}\label{dphterm12}
& \eta(\log x) \int_0^{\eta} x^{-u} H_1(u) \ud u \leq \exp(-2\log\log x) = \frac{1}{(\log x)^2}, \text{ a.s.}
\end{align}
Since we have chosen $\delta\in (0,d_1)$, $0<(1+\delta)/(1+d_1)<1$, and it follows that as $x\to +\infty$,
\begin{align*}
& \eta \log x  = \exp\left\{-\left[D_1^{1/(1+d_1)}(2\log\log x)^{(1+\delta)/(1+d_1)}\right]\right\}\cdot \log x   \\
& = \exp\left[\log \log x - D_1^{\frac{1}{1+d_1}}2^{\frac{1+\delta}{1+d_1}}(\log\log x)^{\frac{1+\delta}{1+d_1}}\right] \\
&\succ \exp\left(\frac{1}{2} \log\log x\right) = \sqrt{\log x} \succ 2\log \log x + \log 2.
\end{align*}
Therefore, we can find a large constant $x_2>0$, such that for all $x\geq x_2$, $\eta \log x > 2\log\log x + \log 2 > 2\log\log x + \eta \log 2$, which implies that
\begin{align}\label{dphterm2}
2(x/2)^{-\eta}& = 2\exp\left(-\eta \log x + \eta\log 2\right) \leq 2\exp\left(-2\log \log x\right) = \frac{2}{(\log x)^2}, \text{ a.s.}
\end{align}
Now we set $x_0=\max(e,x_1,x_2)$. We combine \eqref{boundbypart}, \eqref{dphterm12}, and \eqref{dphterm2} and conclude that for all $x\geq x_0$, $\int_0^{\infty} x^{-u} \ud H_1(u) \leq 3/(\log x)^2$ almost surely. Therefore from \eqref{hfbound1}, we can obtain that for all $x\geq x_0$,
\begin{align*}
\left|h_F(x)\right|& \leq \frac{3}{\underline w_n(\log x)^2}, \text{ a.s.}
\end{align*}
Therefore, Condition (ii) of Theorem \ref{utest1} is satisfied with $B_n=3/\underline w_n$ and $\tau_n=1$. Since Condition (iii) of Theorem \ref{parthm} says that $\overline\alpha_n /\log n \prec \underline w_n \preceq 1$, we have that $B_n\prec \log n/\overline \alpha_n$ and $B_n\prec \log n= \tau_n\log n$. So Condition (iv) of Theorem \ref{utest1} is also verified.

\vspace{.3cm}

\noindent (b) Suppose that $\Pi$ is $\NGGP(a,\kappa,\tau,H_0)$ and the base measure $H_0$ satisfies the conditions in Condition (ii)(b). Then the L\'evy intensity of this NGGP is $\rho(\ud v)=\frac{a}{\Gamma(1-\kappa)}v^{-\kappa-1}e^{-\tau v} \ud v$. From the proof of Theorem \ref{nggpcounter}, we have shown that for sufficiently small $\epsilon_3\in (0,1/2)$, if $h(x)\geq 0$ for $x\in (0,\epsilon_3]$ and $\int_0^{\epsilon} \frac{1}{h(x)^{\kappa}} \ud x < +\infty$, then $\int_0^{\epsilon} \rho[h(x),+\infty) \ud x <+\infty$. Here we take $h(x) = x^{\frac{1}{\kappa(1+\delta)}}$ for some constant $\delta\in (0,(1+d_2)^{1/3}-1)$, where $d_2$ is from Condition (ii)(b). $h(x)$ is convex and increasing in $(0,\epsilon_3]$, and it satisfies $\lim_{x\to 0+}h(x)=0$. Now by exactly the same argument as in the proof for the Dirichlet process in Part (a), we have that $\limsup_{u\to 0+} H_1(u)/h(H_0(u))=0$ almost surely. Hence by Condition (ii)(b), there exists a small constant $0<\epsilon_4 \leq \min(c_2,(\epsilon_3/D_2)^{1/(1+d_2)})$, such that $H_1(u)\leq h(H_0(u))$ almost surely for all $u \in (0,\epsilon_4]$, and meanwhile $H_0(u)\leq D_2 u^{1+d_2}\leq D_2 \epsilon_4^{1+d_2}\leq \epsilon_3$ for all $u \in (0,\epsilon_4]$. Using \eqref{term10}, we have that for all $\eta\in (0,\epsilon_4]$ and any $x\geq e$,
\begin{align}\label{nggpterm11}
& (\log x) \int_0^{\eta} x^{-u} H_1(u) \ud u \leq h\left(D_2\eta^{1+d_2}\right)  =D_2^{1/[\kappa(1+\delta)]} \eta^{(1+d_2)/[\kappa(1+\delta)]} \text{ a.s.}
\end{align}
We choose $\eta = D_2^{-1/(1+d_2)}(\log x)^{-\frac{\kappa(1+\delta)^2}{1+d_2}}$ such that $\eta\leq \epsilon_4$ for all $x\geq x_3$. This holds with the constant $x_3 = \exp\left\{D_2^{-1/[\kappa(1+\delta)^2]}\epsilon_4^{-(1+d_2)/[\kappa(1+\delta)^2]}\right\}$. With this choice of $\eta$, for all $x\geq \max(e,x_3)$, \eqref{nggpterm11} implies that
\begin{align}\label{nggpterm12}
& \eta (\log x) \int_0^{\eta} x^{-u} H_1(u) \ud u \leq 1\cdot D_2^{1/[\kappa(1+\delta)]} \cdot  \left[D_2^{-1/(1+d_2)}(\log x)^{-\frac{\kappa(1+\delta)^2}{1+d_2}}\right]^{(1+d_2)/[\kappa(1+\delta)]} \nonumber \\
&= \frac{1}{(\log x)^{1+\delta}} \text{ a.s.}
\end{align}
Since we have chosen $\delta \in (0,(1+d_2)^{1/3}-1)$, it follows that $0<(1+\delta)^3<1+d_2$. Therefore, as $x\to +\infty$,
\begin{align*}
& \eta \log x  = D_2^{-1/(1+d_2)}(\log x)^{-\frac{\kappa(1+\delta)^2}{1+d_2}}\cdot (\log x) \\
& \geq  D_2^{-1/(1+d_2)} (\log x)^{\frac{(1-\kappa+\delta)(1+\delta)^2}{1+d_2}} \succ (1+\delta)\log\log x + \log 2,
\end{align*}
where the last relation follows since $\kappa\in (0,1)$ and $(1-\kappa+\delta)(1+\delta)^2/(1+d_2)>0$. Therefore, there exists a large constant $x_4\geq e$, such that for all $x\geq x_4$,
\begin{align}\label{nggpterm2}
2(x/2)^{-\eta}& = 2\exp\left(-\eta \log x +\eta \log 2\right) \leq 2\exp\left(-(1+\delta) \log \log x\right) = \frac{2}{(\log x)^{1+\delta}}, \text{ a.s.}
\end{align}
Finally we set $x_0=\max(e,x_3,x_4)$ and combine \eqref{boundbypart}, \eqref{nggpterm12} and \eqref{nggpterm2} to conclude that for all $x\geq x_0$, $\int_0^{\infty} x^{-u} \ud H_1(u) \leq 3/(\log x)^{1+\delta} $ almost surely. Therefore, from \eqref{hfbound1}, we can obtain that for all $x\geq x_0$,
\begin{align*}
\left|h_F(x)\right|& \leq \frac{3}{\underline w_n(\log x)^{1+\delta}}, \text{ a.s.}
\end{align*}
Therefore, Condition (ii) of Theorem \ref{utest1} is satisfied with $B_n=3/\underline w_n$ and $\tau_n=\delta$. Condition (iv) of Theorem \ref{utest1} holds similar to the argument for the DP case in Part (a).
\hfill $\blacksquare$


\section*{Acknowledgements}
We thank Professor Jayanta Ghosh for a comment on the challenges of Bayesian asymptotics for heavy-tailed densities, which served as inspiration for this work. We are grateful to the referees, Associate Editor, and Editor for their comments and suggestions.  LL would like to thank Dong Quan Nguyen for useful discussions on the measurability of the tail index neighborhood. This work was partially supported by National University of Singapore start-up grant R155000172133, NSF grants IIS 1663870 and DMS CAREER 1654579.

\bibliographystyle{Chicago}
\bibliography{tailindex}

\end{document}